\documentclass[11pt]{article}
\usepackage{amssymb,amsmath}
\usepackage{epsfig}
\usepackage{times,bm,mathrsfs}  

\newtheorem{theorem}{Theorem}[section]
\newtheorem{proposition}[theorem]{Proposition}
\newtheorem{definition}[theorem]{Definition}

\newtheorem{lemma}[theorem]{Lemma}

\newtheorem{corollary}[theorem]{Corollary}

\newenvironment{remark}[1][] {
  \smallbreak \noindent {\bf Remark#1:~}} {
  \par\medbreak
}

\newcommand{\qedsymb}{\hfill{\rule{2mm}{2mm}}}
\newenvironment{proof}[1][]{\begin{trivlist}
\item[\hspace{\labelsep}{\bf\noindent Proof#1:\/}]
}{\qedsymb\end{trivlist}}

\def\calP{{\cal P}}

\def\calM{{\cal M}}

\newcommand\ip[1]{{\langle {#1} \rangle}}

\newcommand{\supp}{{\mathsf{supp}}}

\newcommand{\Path}{\mathrm{Path}}

\newcommand{\Star}{\mathrm{Star}}
\newcommand{\Dicttr}{{\cal D}}

\newcommand{\onote}[1]{}
\newcommand{\bnote}[1]{}
\newcommand{\rnote}[1]{}
\newcommand{\enote}[1]{}
\newcommand{\jnote}[1]{}

\newcommand{\eps}{\varepsilon}
\renewcommand{\epsilon}{\varepsilon}

\newcommand{\ignore}[1]{}

\def \bits    {\{-1,1\}}

\def \bn      {\bits^n}
\def \bm      {\bits^m}

\def \E       {{\bf E}}

\def \eps     {\epsilon}

\def \half    {{\textstyle \frac{1}{2}}}

\def \infint  {\int_{-\infty}^\infty}

\def \la      {\langle}
\def \ra      {\rangle}
\def \maj     {\mathsf{MAJ}}

\def \P       {{\bf P}}

\renewcommand{\Pr}{\P}

\def \Reals   {\mathbb{R}}

\begin{document}

\title{\bf Non-interactive correlation distillation, inhomogeneous Markov
chains, and the reverse Bonami-Beckner inequality}

\author{
 Elchanan Mossel \footnote{Department of Statistics, U.C. Berkeley. {\tt mossel@stat.berkeley.edu}. Supported by a
 Miller fellowship in CS and Statistics, U.C. Berkeley.}
\and
 Ryan O'Donnell \footnote{Institute for Advanced Study, Princeton, NJ.  {\tt
   odonnell@ias.edu}. Most of this work was done while the author was a student at
   Massachusetts Institute of Technology.  This material is based upon work supported by the
   National Science Foundation under agreement No.\ CCR-0324906.  Any opinions, findings and
   conclusions or recommendations expressed in this material are those of the authors and do
   not necessarily reflect the views of the National Science Foundation.}
 \and
 Oded Regev \footnote{Department of Computer Science, Tel-Aviv University, Tel-Aviv 69978, Israel.
   Most of this work was done while the author was at the
   Institute for Advanced Study, Princeton, NJ.
   Work supported by an Alon Fellowship, ARO grant DAAD19-03-1-0082 and
   NSF grant CCR-9987845.}
 \and
 Jeffrey E. Steif \footnote{Department of Mathematics, Chalmers University of Technology,
 412 96 Gothenburg, Sweden. {\tt steif@math.chalmers.se}.
Supported in part by NSF grant DMS-0103841 and in part by
the Swedish Research Council.}
 \and
 Benny Sudakov \footnote{Department of Mathematics,
Princeton University, Princeton, NJ 08544, USA.
{\tt bsudakov@math.princeton.edu}. Research supported in part by NSF
grant DMS-0106589, DMS-0355497, and by an Alfred P. Sloan fellowship.}
}

\date{}

\maketitle

\begin{abstract}

In this paper we study \emph{non-interactive
correlation distillation (NICD)}, a generalization of noise
sensitivity previously considered in~\cite{AMW91, MO03, Yang04}.
We extend the model to \emph{NICD on trees}.  In this model there
is a fixed undirected tree with players at some of the nodes.  One
node is given a uniformly random string and this string is
distributed throughout the network, with the edges of the tree
acting as independent binary symmetric channels.  The goal of the
players is to agree on a shared random bit without communicating.

Our new contributions include the following:

\begin{itemize}

\item In the case of a $k$-leaf star graph (the model considered
in \cite{MO03}), we resolve the open question of whether the
success probability must go to zero as $k \to \infty$.  We show
that this is indeed the case and provide matching upper and lower
bounds on the asymptotically optimal rate (a slowly-decaying
polynomial).

\item In the case of the $k$-vertex path graph, we show
that it is always optimal for all players to use the same
1-bit function.

\item In the general case we show that all players should use
monotone functions. We also show, somewhat surprisingly, that
for certain trees it is better if not all players use the same function.

\end{itemize}
Our techniques include the use of the \emph{reverse} Bonami-Beckner
inequality.  Although the usual Bonami-Beckner has been frequently used
before,
its reverse counterpart seems very little-known;
To demonstrate its strength, we use
it to prove a new isoperimetric inequality for the discrete cube
and a new result on the mixing of short random walks on the cube.
Another tool that we need is a tight bound on the
probability that a Markov chain stays inside certain sets; we
prove a new theorem generalizing and strengthening previous such
bounds~\cite{AjtaiKS, AlonFWZ,AlonSpencer}.
On the probabilistic side, we use the ``reflection principle'' and the FKG
and related inequalities in order to study the problem on general trees.
\end{abstract}

\thispagestyle{empty}
\newpage
\addtocounter{page}{-1}

\section{Introduction}

\subsection{Non-interactive correlation --- the problem and
previous work} Our main topic in this paper is the problem of
\emph{non-interactive correlation distillation (NICD)}, previously
considered in~\cite{AMW91, MO03, Yang04}.  In its most general
form the problem involves $k$ players who receive noisy copies of
a uniformly random bit string of length $n$.  The players wish to
agree on a single random bit but are not allowed to communicate.
The problem is to understand the extent to which the players can
successfully distil the correlations in their strings into a
shared random bit.  This problem is relevant for cryptographic
information reconciliation, random beacons in cryptography and
security, and coding theory; see~\cite{Yang04}.

In its most basic form, the problem involves only two players; the
first gets a uniformly random string $x$ and the second gets a
copy $y$ in which each bit of $x$ is flipped independently with
probability $\eps$.  If the players try to agree on a shared bit
by applying the same Boolean function $f$ to their strings, they
will fail with probability $\Pr[f(x) \neq f(y)]$.  This quantity
is known as the \emph{noise sensitivity of $f$ at $\epsilon$}, and
the study of noise sensitivity has played an important role in
several areas of mathematics and computer science (e.g.,
inapproximability~\cite{Has01}, learning
theory~\cite{BJT99a,KOS02}, hardness amplification~\cite{Odo02},
mixing of short random walks~\cite{KKL88},
percolation~\cite{BKS99}; see also~\cite{Odo03}). In~\cite{AMW91},
Alon, Maurer, and Wigderson showed that if the players want to use
a balanced function $f$, no improvement over the naive strategy of
letting $f(x) = x_1$ can be achieved.

The paper~\cite{MO03} generalized from the two-player problem NICD
to a $k$-player problem, in which a uniformly random string $x$ of
length $n$ is chosen, $k$ players receive independent
$\eps$-corrupted copies, and they apply (possibly different)
balanced Boolean functions to their strings, hoping that all
output bits agree.  This generalization is equivalent to studying
high norms of the Bonami-Beckner operator applied to Boolean
functions (i.e., $\|T_\rho f\|_k$); see Section~\ref{sec:bb} for
definitions.  The results in~\cite{MO03} include:
optimal protocols involve all players using the same function;
optimal functions are always monotone; for $k = 3$ the first-bit
(`dictator') is best; for fixed $\eps$ and fixed $n$ and $k \to \infty$,
all players should use the majority function; and, for fixed $n$
and $k$ and $\eps \to 0$ or $\eps \to 1/2$ dictator is best.

Later, Yang~\cite{Yang04} considered a different generalization of
NICD, in which there are only two players but the corruption model
is different from the ``binary symmetric channel'' noise
considered previously.  Yang showed that for certain more general
noise models, it is still the case that the dictator function is
optimal; he also showed an upper bound on the players' success rate in
the erasure model.

\subsection{NICD on trees; our results}
In this paper we propose a natural generalization of the NICD
models of~\cite{AMW91, MO03}, extending to a tree topology.  In
our generalization we have a network in the form of a tree; $k$ of
the nodes have a `player' located on them. One node broadcasts a
truly random string of length $n$. The string follows the edges of
the trees and eventually reaches all the nodes. Each edge of the
tree independently introduces some noise, acting as a binary
symmetric channel with some fixed crossover probability $\eps$.
Upon receiving their strings, each player applies a balanced
Boolean function, producing one output bit.  As usual, the goal of
the players is to agree on a shared random bit without any further
communication; the protocol is successful if all $k$ parties
output the same bit. (For formal definitions, see
Section~\ref{sec:prelims}.)  Note that the problem considered
in~\cite{MO03} is just NICD on the star graph of $k+1$ nodes with
the players at the $k$ leaves.\\

We now describe our new results:

\paragraph{The $k$-leaf star graph:} We first study the same $k$-player
star problem considered in~\cite{MO03}.  Although this paper found
maximizing protocols in certain asymptotic scenarios for the
parameters $k$, $n$, and $\eps$, the authors left open
what is arguably the most interesting setting: $\eps$ fixed,
$k$ growing arbitrarily large, and $n$ unbounded in terms of
$\eps$ and $k$.  Although it is natural to guess that the success
rate of the players must go to zero exponentially fast in terms of
$k$, this turns out not to be the case; \cite{MO03} notes that if
all players apply the majority function (with $n$ large enough)
then they succeed with probability $\Omega(k^{-C(\eps)})$ for some
finite constant $C(\eps)$ (the estimate~\cite{MO03} provides is
not sharp).  \cite{MO03} left as a major
open problem to prove that the success probability goes to $0$ as
$k \to \infty$.

In this paper we solve this problem.  In Theorem \ref{thm:exactrate} we
show that the success probability must indeed go to zero as $k \to
\infty$. Our upper bound is a slowly-decaying polynomial.
Moreover, we provide a matching lower bound: this follows from a
tight analysis of the majority protocol. The proof of our upper
bound depends crucially on the reverse Bonami-Beckner inequality,
an important tool that will be described later.

\paragraph{The $k$-vertex path graph:}
In the case of NICD on the path graph, we prove in Theorem
\ref{thm:const_prot} that in the optimal protocol
all players should use the same
1-bit function. In order to prove this, we prove in Theorem \ref{thm:aks},
a new tight bound
on the probability that a Markov chain stays inside certain sets.
Our theorem generalizes and strengthens previous
work~\cite{AjtaiKS, AlonFWZ, AlonSpencer}.

\paragraph{Arbitrary trees:} In this general case, we show in Theorem
\ref{thm:optimal_monotone} that
there always exists an optimal protocol in which all players use
monotone functions. Our analysis uses methods of discrete
symmetrization together with the FKG correlation inequality.

In Theorem \ref{thm:no_const_prot} we show that for certain trees it is better
if not all players use the same function. This might be
somewhat surprising: after all, if all players wish
to obtain the same result, won't they be better off using
the same function? The intuitive reason the answer to this
is negative can be explained by Figure \ref{fig:starpath}: players on
the path and players on the star each `wish' to use a
different function. Those on the star wish to use the majority
function and those on the path wish to use a dictator function.
Indeed, we will show that this strategy yields better
success probability than any strategy in which all players use the
same function.

\subsection{The reverse Bonami-Beckner inequality}
We would like to
highlight the use of the reverse Bonami-Beckner inequality,
mentioned above.
Let us start by describing the original Bonami-Beckner inequality,
see Theorem \ref{thm:bonami-beckner}.
This inequality considers an operator known as the Bonami-Beckner
operator (see Section~\ref{sec:bb}). It says that some high norm
of the result of the Bonami-Beckner operator applied to a function
can be upper-bounded by some low norm of the original function.
Its main strength is in its ability to relate two different norms;
this is the reason it is often referred to as a hypercontractive
inequality. The inequality was originally proved by Bonami in 1970
\cite{Bon70} and then independently by Beckner in 1973
\cite{Bec75}. It was first used to analyze discrete problems in a
a remarkable paper by Kahn, Kalai and Linial
\cite{KKL88} where they considered the influence of variables on
Boolean functions. The inequality has proved to be of great
importance in the study of combinatorics of
$\{0,1\}^n$~\cite{BKKKL92,BK97b,Fri98},
percolation and random graphs~\cite{Tal94,FK96,BKS99,Bou99} and many
other applications~\cite{BL90,AKRS94,RazFourier,AmanoMaruoka,ODonnellServedio,
DGK02,DS02,Kho02,Odo02}.


Far less well-known is the fact that the Bonami-Beckner inequality
admits a reversed form.
This reversed form
was first proved by Borell~\cite{Bor82} in 1982. Unlike the
original inequality, the reverse inequality says that some low
norm of the Bonami-Beckner operator applied to a non-negative
function can be bounded {\em below} by some higher norm of the
original function. Moreover, the norms involved in the reverse
inequality are all at most $1$ while the norms in the original
inequality are all at least $1$. A final difference between the
two inequalities is that in the reverse inequality we need to
assume that the function is non-negative.

We are not aware of any previous uses of the reverse Bonami-Beckner inequality
for the study of discrete problems. The inequality
seems very promising and we hope it will prove useful in the future. To demonstrate its
strength, we provide two applications:

\paragraph{Isoperimetric inequality on the discrete cube:} As a corollary of the reverse
Bonami-Beckner inequality, we obtain
in Theorem \ref{thm:isop} an isoperimetric inequality
on the discrete cube. Although it is a simple corollary, we
believe that the isoperimetric inequality is interesting. It is
also used later to give a sort of hitting time upper-bound for
short random walks. In order to illustrate it, let us consider two
subsets $S, T\subseteq \bn$ each containing a constant fraction $\sigma$ of
the $2^n$ elements of the discrete cube. We now perform the
following experiment: we choose a random element of $S$ and flip
each of its $n$ coordinates with probability $\eps$ for some small
$\eps$. What is the probability that the resulting element is in
$T$? Our isoperimetric inequality implies that it is at least some
constant independent of $n$.
For example, given any two sets with fractional size $1/3$, the
probability that flipping each coordinate with probability $.3$
takes a random point chosen from the first set into the second set
is at least $(1/3)^{1.4/.6} \approx 7.7\%$. We also show that our
bound is close to tight. Namely, we analyze the above
probability for diametrically opposed Hamming balls and show that
it is close to our lower bound.

\paragraph{Short random walks:} Our second application in Proposition
\ref{prop:randomwalks} is to
short random walks on the
discrete cube. Consider the following scenario. We have two sets
$S,T \subseteq \bn$ of size at least $\sigma 2^n$ each. We start a
walk from a random element of the set $S$ and at each time step
proceed with probability $1/2$ to one of its neighbors which we
pick randomly. Let $\tau n$ be the length of the random walk. What
is the probability that the random walk terminates in $T$? If
$\tau = C \log n$ for a large enough constant $C$ then it is known
that the random walk mixes and therefore we are guaranteed to be
in $T$ with probability roughly $\sigma$. However, what happens if
$\tau$ is, say, $0.2$? Notice that $\tau n$ is then less than the
diameter of the cube! For certain sets $S$, the random walk might
have zero probability to reach certain vertices, but if $\sigma$
is at least, say, a constant then there will be some nonzero
probability of ending in $T$.  We bound from below the probability
that the walk ends in $T$ by a function of $\sigma$ and $\tau$
only. For example, for $\tau = 0.2$, we obtain a bound of roughly
$\sigma^{10}$. The proof crucially depends on the reverse
Bonami-Beckner inequality; to the best of our knowledge, known
techniques, such as spectral methods, cannot yield a similar
bound.

\section{Preliminaries} \label{sec:prelims}
We now formally define the problem of ``non-interactive
correlation distillation (NICD) on trees with the binary symmetric
channel (BSC).''  In general we have four parameters.  The first
is $T$, an undirected tree giving the geometry of the problem.
Later the vertices of $T$ will become labeled by binary
strings, and the edges of $T$ will be thought of as independent
binary symmetric channels. The second parameter of the problem is
$0 < \rho < 1$ which gives the \emph{correlation} of bits on opposite
sides of a channel.  By this we mean that if a bit string $x \in
\bn$ passes through the channel producing the bit string $y \in
\bn$ then $\E[x_i y_i] = \rho$ independently for each $i$. We say
that $y$ is a $\rho$-correlated copy of $x$.  We will also
sometimes refer to $\eps = \half - \half \rho \in (0,\half)$,
which is the probability with which a bit gets flipped
--- i.e., the crossover probability of the channel. The third
parameter of the problem is $n$, the number of bits in the string
at every vertex of $T$.  The fourth parameter of the problem is a
subset of the vertex set of $T$, which we denote by $S$. We refer
to the $S$ as the set of \emph{players}. Frequently $S$ is simply
all of $V(T)$, the vertices of $T$.

To summarize, an instance of the NICD on trees problem is
parameterized by:

\begin{enumerate}
\item[1.]
$T$, an undirected tree;
\item[2.] $\rho \in
(0,1)$, the correlation parameter;
\item[3.] $n \geq 1$, the
string length; and,
\item[4.] $S \subseteq V(T)$, the set of players.
\end{enumerate}

Given an instance, the following process happens.  Some vertex $u$
of $T$ is given a uniformly random string $x^{(u)} \in \bn$.  Then
this string is passed through the BSC edges of $T$ so that every
vertex of $T$ becomes labeled by a random string in $\bn$.  It is
easy to see that the choice of $u$ does not matter, in the sense
that the resulting joint probability distribution on strings for
all vertices is the same regardless of $u$.  Formally speaking, we
have $n$ independent copies of a ``tree-indexed Markov chain;'' or
a ``Markov chain on a tree'' \cite{Georgii:88}.
The index set is $V(T)$ and the probability measure $\P$ on
$\alpha \in \bits^{V(T)}$ is defined by
\[
\P(\alpha) = \half \left(\half + \half \rho\right)^{A(\alpha)}
\left(\half - \half \rho\right)^{B(\alpha)},
\]
where $A(\alpha)$ is the number of pairs of neighbors where
$\alpha$ agrees and $B(\alpha)$ is the number of pairs of
neighbors where $\alpha$ disagrees.

Once the strings are distributed on the vertices of $T$, the
player at the vertex $v \in S$ looks at the string $x^{(v)}$ and
applies a (pre-selected) Boolean function $f_v \colon \bn \to
\bits$.  The goal of the players is to maximize the probability
that the bits $f_v(x^{(v)})$ are identical for \emph{all} $v \in
S$.  In order to rule out the trivial solutions of constant
functions and to model the problem of flipping a shared random
coin, we insist that all functions $f_v$ be \emph{balanced}; i.e.,
have equal probability of being $-1$ or $1$.  As noted in
\cite{MO03}, this does not necessarily ensure that when all
players agree on a bit it is conditionally equally likely to be
$-1$ or $1$; however, if the functions are in addition
antisymmetric, this property does hold.  We call a collection of
balanced functions $(f_v)_{v \in S}$ a \emph{protocol} for the
players $S$, and we call this protocol \emph{simple} if all of the
functions are the same.

To conclude our notation, we write $\calP(T, \rho, n, S, (f_v)_{v \in S})$
for the probability that the protocol succeeds -- i.e., that all
players output the same bit.  When the protocol is simple we write
merely $\calP(T, \rho, n, S, f)$.  Our goal is to study the maximum this
probability can be over all choices of protocols. We denote by
\[
\calM(T, \rho, n, S) = \sup_{(f_v)_{v \in S}} \calP(T, \rho, n, S,
(f_v)_{v \in S}),
\]
and define
\[
\calM(T, \rho, S) = \sup_n \calM(T, \rho, n, S).
\]

\section{Reverse Bonami-Beckner and applications} \label{sec:bb}

In this section we recall the little-known reverse Bonami-Beckner
inequality and obtain as a corollary
an isoperimetric inequality on the discrete cube.  These results will
be useful in analyzing the NICD problem on the star graph and we
believe they are of independent interest. We also
obtain a new result about the mixing of relatively short random
walks on the discrete cube.

\subsection{The reverse Bonami-Beckner inequality} \label{sec:wherearewe}
We begin with a discussion of the Bonami-Beckner inequality.
Recall the Bonami-Beckner operator $T_\rho$, a linear operator on
the space of functions $\bn \to \Reals$ defined by
\[
(T_\rho f)(x) = \E[f(y)],
\]
where $y$ is a $\rho$-correlated copy of $x$.  The usual
Bonami-Beckner inequality, first proved by Bonami~\cite{Bon70} and
later independently by Beckner~\cite{Bec75}, is the following:

\begin{theorem} \label{thm:bonami-beckner}
Let $f \colon \bn \to \Reals$ and $q \geq p \geq 1$.  Then
\[
\|T_\rho f\|_q \leq \|f\|_p \qquad \text{for all $0 \leq \rho
\leq (p-1)^{1/2}/(q-1)^{1/2}$.}
\]
\end{theorem}

The reverse Bonami-Beckner inequality is the following:

\begin{theorem} \label{thm:reverse-BB}
Let $f \colon \bn \to \Reals^{\geq 0}$ be a nonnegative function
and let $-\infty < q \leq p \leq 1$. Then
\begin{equation} \label{eq:borell}
\|T_\rho f\|_q \geq \|f\|_p \qquad \text{for all $0 \leq \rho
\leq (1-p)^{1/2}/(1-q)^{1/2}$.}
\end{equation}
\end{theorem}

Note that in this theorem we consider $r$-norms for $r \leq 1$.
The case of $r = 0$ is a removable singularity: by $\|f\|_0$ we
mean the geometric mean of $f$.  Note also that since $T_\rho$ is
a convolution operator, it is positivity-improving for any $\rho <
1$; i.e., when $f$ is nonnegative so too is $T_\rho f$, and if $f$
is further not identically zero, then $T_\rho f$ is everywhere
positive.

The reverse Bonami-Beckner theorem is proved in the same way the
usual Bonami-Beckner theorem is proved; namely, one proves the
inequality in the case of $n=1$ by elementary means, and then
observes that the inequality tensors. Since Borell's original proof may be too
compact to be read by some, we provide an expanded version of it
in Appendix~\ref{sec:proof_rev_beckner} for completeness.

We will actually need the following ``two-function'' version of
the reverse Bonami-Beckner inequality which follows easily from
the reverse Bonami-Beckner inequality using the (reverse)
H\"{o}lder inequality (see Appendix~\ref{sec:proof_rev_beckner}):

\begin{corollary} \label{cor:two-fcn-reverse-BB} Let $f, g \colon \bn \to \Reals^{\geq 0}$ be nonnegative, let
$x \in \bn$ be chosen uniformly at random, and let $y$ be a
$\rho$-correlated copy of $x$.  Then for $-\infty < p,q < 1$,
\begin{equation} \label{eq:two_Borell}
\E[f(x)g(y)] \geq \|f\|_p \|g\|_q \qquad \text{for all $0 \leq
\rho \leq (1-p)^{1/2}(1-q)^{1/2}$.}
\end{equation}
\end{corollary}

\subsection{A new isoperimetric inequality on the discrete cube}

In this subsection we use the reverse Bonami-Beckner inequality to
prove an isoperimetric inequality on the discrete
cube.  Let $S$ and $T$ be two subsets of $\bn$. Suppose that $x
\in \bn$ is chosen uniformly at random and $y$ is a
$\rho$-correlated copy of $x$. We obtain the following theorem,
which gives a lower bound on the probability that $x \in S$ and $y
\in T$ as a function of $|S|/2^n$ and $|T|/2^n$ only.

\begin{theorem} \label{thm:isop}
Let $S, T \subseteq \bn$ with $|S| = \exp(-s^2/2)2^n$ and $|T| =
\exp(-t^2/2)2^n$.  Let $x$ be chosen uniformly at random from
$\bn$ and let $y$ be a $\rho$-correlated copy of $x$.  Then
\begin{equation} \label{eq:isop}
\Pr[x \in S, y \in T] \geq
\exp\left(-\frac12 \frac{s^2 + 2\rho s t + t^2}{1 - \rho^2}\right).
\end{equation}
\end{theorem}
\begin{proof}
Take $f$ and $g$ to be the $0$-$1$ characteristic functions of $S$
and $T$, respectively.  Then by
Corollary~\ref{cor:two-fcn-reverse-BB}, for any choice of $p,
q < 1$ with $(1-p)(1-q) = \rho^2$, we get
\begin{equation}
\Pr[x \in S, y \in T] = \E[f(x)g(y)] \geq \|f\|_p \|g\|_q = \exp(-s^2/2p) \exp(-t^2/2q).
\label{eqn:iso}
\end{equation}
Write $p = 1 - \rho r$, $q = 1 - \rho/r$ in~(\ref{eqn:iso}), with
$r > 0$.  Maximizing the right-hand side as a function of $r$
the best choice is $r = ((t/s) +
\rho)/(1 + \rho(t/s))$ which yields in turn
\[
\begin{array}{ll}
p = 1-\rho r = \frac{1-\rho^2}{1+\rho (t/s)}, &
q = 1-\rho/r = \frac{t}{s} \frac{1-\rho^2}{\rho + (t/s)}.
\end{array}
\]
(Note that this depends only on the ratio
of $t$ and $s$.)  Substituting this choice of $r$ (and hence $p$
and $q$) into~(\ref{eqn:iso}) yields $\exp(-\half \frac{s^2 +
2\rho st + t^2}{1-\rho^2})$, as claimed.
\end{proof}

\vspace{0.15cm}

Taking $\sigma = \exp(-s^2/2)$ and $\sigma^\alpha = \exp(-t^2/2)$
we obtain
\[
-\frac12(s^2 + 2\rho s t + t^2)=
\log \sigma - \rho \sqrt{-2 \log \sigma} \sqrt{-2 \alpha \log \sigma} +
\alpha \log \sigma = \log \sigma (1 + 2 \rho \sqrt{\alpha} + \alpha),
\]
and therefore
\[
\exp\left(-\frac12 \frac{s^2 + 2\rho s t + t^2}{1 - \rho^2}\right) =
\sigma^{(1 + 2 \rho \sqrt{\alpha} + \alpha)/(1-\rho^2)}.
\]
Thus, conditioned on starting at $S$, the probability of ending at $T$ is
at least
\[
\sigma^{(1 + 2 \rho \sqrt{\alpha} + \alpha)/(1-\rho^2)-1}=
\sigma^{(\sqrt{\alpha} + \rho)^2/(1-\rho^2)}.
\]
We thus obtain the following corollary of Theorem~\ref{thm:isop}.
\begin{corollary} \label{cor:isop}
Let $S \subseteq \bn$ have fractional size $\sigma \in [0,1]$, and
let $T \subseteq \bn$ have fractional size $\sigma^\alpha$, for
$\alpha \geq 0$.  If $x$ is chosen uniformly at random from $S$
and $y$ is a $\rho$-correlated copy of $x$, then the probability
that $y$ is in $T$ is at least
\[
\sigma^{(\sqrt{\alpha} + \rho)^2/(1-\rho^2)}.
\]
In particular, if $|S| = |T|$ then this probability is at least
$\sigma^{(1+\rho)/(1-\rho)}$.
\end{corollary}

In Subsection \ref{subsec:tightness} below we show that the isoperimetric
inequality is almost tight. First, we prove a similar bound for random
walks on the cube.


\subsection{Short random walks on the discrete cube}

We can also prove a result of a similar flavor about short random walks on the discrete cube:

\begin{proposition}\label{prop:randomwalks}
Let $\tau>0$ be arbitrary and let $S$ and $T$ be two subsets of $\bn$.
Let $\sigma \in [0,1]$ be the fractional size of $S$ and let $\alpha$
be such that the fractional size of $T$ is $\sigma^\alpha$.
Consider a standard random walk on the discrete cube that starts from
a uniformly random vertex in $S$ and walks for $\tau n$ steps.
Here by a standard random
walk we mean that at each time step we do nothing
with probability $1/2$ and we walk along the $i$th edge with
probability $1/2n$. Let $p^{(\tau n)}(S, T)$ denote the probability that the walk
ends in $T$. Then,
$$p^{(\tau n)}(S, T) \geq
   \sigma^{\frac{(\sqrt{\alpha} + \exp(-\tau))^2}{1-\exp(-2\tau)}} -
O\Big(\frac{\sigma^{(-1 + \alpha)/2}}{\tau n}\Big).
$$
In particular, when $|S| = |T| = \sigma 2^n$ then $
p^{(\tau n)}(S,T) \geq
\sigma^{\frac{1+\exp(-\tau)}{1-\exp(-\tau)}} - O(\frac{1}{\tau n})$.
\end{proposition}

\noindent
The Laurent series of $\frac{1+e^{-\tau}}{1-e^{-\tau}}$ is $2/\tau + \tau/6 - O(\tau^3)$ so for $1/\log n \ll \tau \ll 1$ our bound is roughly $\sigma^{2/\tau}$.

For the proof we will first need a simple lemma:
\begin{lemma}\label{clm:easytosee}
For $y>0$ and any $0 \le x \le y$,
$$ 0 \le e^{-x} - (1-x/y)^y \le O(1/y).$$
\end{lemma}
\begin{proof}
The expression above can be written as
$$ e^{-x} - e^{y \log (1-x/y)}.$$
We have $\log(1-x/y) \le -x/y$ and hence we obtain the first inequality.
For the second inequality, notice that if $x \ge 0.1y$ then
both expressions are of the form $e^{-\Omega(y)}$ which is certainly $O(1/y)$.
On the other hand, if $0\le x < 0.1y$ then there is a constant $c$ such that
$$\log(1-x/y) \ge -x/y - cx^2/y^2.$$
The Mean Value Theorem implies that for $0\le a \le b$, $e^{-a}-e^{-b} \le e^{-a}(b-a)$.
Hence,
$$ e^{-x} - e^{y \log (1-x/y)} \le e^{-x}(- y \log (1-x/y) - x) \le \frac{c x^2e^{-x}}{y}.$$
The lemma now follows because $x^2 e^{-x}$ is uniformly bounded
for $x \geq 0$.
\end{proof}

We now prove Proposition~\ref{prop:randomwalks}.  The proof uses Fourier analysis; for the required definitions see, e.g., \cite{KKL88}.

\begin{proof}
Let $x$ be a uniformly random point in $\bn$ and $y$ a point
generated by taking a random walk of length $\tau n$ starting from
$x$.  Let $f$ and $g$ be the $0$-$1$ indicator functions of $S$
and $T$, respectively, and say $E[f] = \sigma$, $E[g] =
\sigma^\alpha$. Then by writing $f$ and $g$ in their Fourier
decomposition we obtain that
$$\sigma \cdot p^{(\tau n)}(S, T)  =  \Pr[x \in S, y \in T] =  \E[f(x)g(y)]
 =  \sum_{U, V} \hat{f}(U) \hat{g}(V) \E[x_U y_V]$$
where $U$ and $V$ range over all subsets of $\{1,\ldots,n\}$. Note that  $\E[x_U y_V]$ is zero unless $U=V$. Therefore
\begin{eqnarray*}
\sigma p^{(\tau n)}(S, T) &= & \sum_U \hat{f}(U) \hat{g}(U) \E[(xy)_U]
 =  \sum_U \hat{f}(U) \hat{g}(U) \Big(1-\frac{|U|}{n}\Big)^{\tau n}  \\
& = & \sum_U \hat{f}(U) \hat{g}(U) \exp(-\tau |U|) + \sum_U
\hat{f}(U) \hat{g}(U) \Big[\Big(1-\frac{|U|}{n}\Big)^{\tau n} - \exp(-\tau |U|)\Big] \\
& = & \la f, T_{\exp(-\tau)} g \ra + \sum_U \hat{f}(U) \hat{g}(U)
\Big[\Big(1-\frac{|U|}{n}\Big)^{\tau n} - \exp(-\tau |U|)\Big] \\
& \geq & \la f, T_{\exp(-\tau)} g \ra  -
\max_{|U|} \Big| \Big(1-\frac{|U|}{n}\Big)^{\tau n} - \exp(-\tau |U|) \Big|
\sum_U |\hat{f}(U)
\hat{g}(U)|.
\end{eqnarray*}
By Corollary~\ref{cor:isop},
\[
\sigma^{-1}\la f, T_{\exp(-\tau)} g \ra \geq
\sigma^{\frac{(\sqrt{\alpha} + \exp(-\tau))^2}{1-\exp(-2\tau)}}.
\]
By Cauchy-Schwarz and Parseval's identity,
$$\sum_U |\hat{f}(U) \hat{g}(U)| \leq \|\hat{f}\|_2
\|\hat{g}\|_2 = \|f\|_2 \|g\|_2 = \sigma^{(1+\alpha)/2}.$$
In addition, from Lemma~\ref{clm:easytosee} with $x=\tau |U|$ and $y=\tau n$ we have
that
\[
\max_{|U|} \Big|\Big(1-\frac{|U|}{n}\Big)^{\tau n} - \exp(-\tau |U|) \Big| =
O\Big(\frac{1}{\tau n}\Big).
\]
Hence,
\[
p^{(\tau n)}(S,T) \geq \sigma^{\frac{(\sqrt{\alpha} +
\exp(-\tau))^2}{1-\exp(-2\tau)}} - O\Big(\frac{\sigma^{(-1 + \alpha)/2}}{\tau
n}\Big).
\]
\end{proof}

\subsection{Tightness of the isoperimetric inequality} \label{subsec:tightness}
We now show that Theorem~\ref{thm:isop} is almost tight.
Suppose $x \in\bn$ is chosen uniformly at random and $y$ is a $\rho$-correlated copy of $x$.
Let us begin by understanding more about how $x$ and $y$
are distributed. Define
\[
\Sigma(\rho) = \left[\begin{array}{cc}1 & \rho
\\ \rho & 1 \end{array}\right]
\]
and recall that the density function of the bivariate normal distribution
$\phi_{\Sigma(\rho)} : \Reals^2 \to \Reals^{\geq 0}$ with mean $0$ and
covariance matrix $\Sigma(\rho)$, is given by
\begin{eqnarray*}
\phi_{\Sigma(\rho)}(x,y) &=&
(2\pi)^{-1}(1-\rho^2)^{-\half}\exp\left(-\frac12\frac{x^2 - 2\rho
x y + y^2}{1-\rho^2}\right) \\
& = & (1-\rho^2)^{-\half} \phi(x)\phi\left(\frac{y-\rho
x}{(1-\rho^2)^\half}\right).
\end{eqnarray*}
Here $\phi$ denotes the standard normal density function on
$\Reals$, $\phi(x)=(2\pi)^{-1/2} e^{-x^2/2}$.\\

\begin{proposition}  \label{prop:normallimit} Let $x \in \bn$ be chosen uniformly at random,
and let $y$ be a $\rho$-correlated copy of $x$.  Let $X =
n^{-1/2}\sum_{i=1}^n x_i$ and $Y = n^{-1/2}\sum_{i=1}^n y_i$. Then
as $n \to \infty$, the pair of random variables $(X,Y)$ approaches
the distribution $\phi_{\Sigma(\rho)}$.
As an error bound, we have that for any convex region $R \subseteq
\Reals^2$,
$$\left|\Pr\big[(X,Y) \in R\big] - \int\!\!\!\int_{R} \phi_{\Sigma(\rho)}(x,y)\,dy\,dx\right| \leq
O((1-\rho^2)^{-1/2}n^{-1/2}).$$
\end{proposition}
\begin{proof}
This follows from the Central Limit Theorem (see,
e.g.,~\cite{Fel68}), noting that for each coordinate $i$,
$\E[x_i^2] = \E[y_i^2] = 1$, $\E[x_i y_i] = \rho$.  The
Berry-Ess\'een-type error bound is proved in Sazonov~\cite[p.\@ 10,
Item 6]{Saz81}.
\end{proof}

Using this proposition we can obtain the following result
for two diametrically opposed Hamming balls.

\begin{proposition} \label{prop:almost}
Fix $s,t > 0$, and let $S,T \subseteq \bn$ be diametrically
opposed Hamming balls, with $S = \{x \colon \sum_i x_i$
$\leq -sn^{1/2} \}$ and $T = \{x \colon \sum_i x_i \geq t n^{1/2} \}$. Let $x$ be chosen uniformly at random from
$\bn$ and let $y$ be a $\rho$-correlated copy of $x$. Then
we have
\[
\lim_{n \to \infty} \Pr[x \in S, y \in T] \leq
\frac{\sqrt{1-\rho^2}}{2 \pi s (\rho s + t)}
\exp\left(-\frac12 \frac{s^2 + 2\rho s t + t^2}{1 - \rho^2}\right).
\]
\end{proposition}

\begin{proof}
%
\begin{eqnarray*}
\lim_{n \to \infty} \Pr[x \in S, y \in T] &=&
\int_{s}^\infty \int_{t}^\infty
\phi_{\Sigma(-\rho)}(x,y)\,dy\,dx\,\,\,
\left(\mbox{ By Lemma \ref{prop:normallimit} } \right) \\
& \leq &
\int_s^{\infty} \int_t^{\infty}
\frac{x(\rho x + y)}{s (\rho s + t)} \phi_{\Sigma(-\rho)}(x,y)\,dy\,dx \\
&\,&
\left(\mbox{  Since } \frac{x(\rho x + y)}{s(\rho s + t)} \geq 1 \mbox{ on }
x \geq s, y  \geq t \right) \\
&=&
\frac{1}{\sqrt{1-\rho^2}}\int_s^{\infty} \int_t^{\infty}
\frac{x(\rho x + y)}{s(\rho s + t)} \phi(x)
\phi\left(\frac{y+\rho x}{\sqrt{1-\rho^2}}\right) \,dy\,dx \\
& \leq &
\frac{1}{\sqrt{1-\rho^2}}
\int_s^{\infty} \int_{\rho s + t}^{\infty} \frac{x z}{s(\rho s + t)} \phi(x)
\phi\left(\frac{z}{\sqrt{1-\rho^2}}\right) \,dz\,dx \\
&\,&
\left( \mbox{ Using } z = \rho x + y \mbox{ and noting }
\frac{x z}{s(\rho s + t)} \geq 1 \mbox{ on }
x \geq s, z \geq \rho s + t \right)
\\ &=&
\frac{1}{s (\rho s + t)\sqrt{1-\rho^2}}\left(\int_s^{\infty} x \phi(x)
  dx
\right)\left(\int_{\rho s + t}^{\infty} z
\phi\left(\frac{z}{\sqrt{1-\rho^2}}\right) \,dz \right)
\\ &=&
\frac{\sqrt{1-\rho^2}}{s (\rho s + t)} \phi(s)
\phi\left(\frac{\rho s + t}{\sqrt{1-\rho^2}}\right)
\\ &=& \frac{\sqrt{1-\rho^2}}{2 \pi s (\rho s + t)}
\exp\left(-\frac12 \frac{s^2 + 2\rho s t + t^2}{1 - \rho^2}\right).
\end{eqnarray*}
The result follows.
\end{proof}



\vspace{0.1cm} \noindent By the Central Limit Theorem, the set $S$
in the above statement satisfies (see~\cite[26.2.12]{AS72}),
$$\lim_{n\to \infty} |S|2^{-n}
=\frac{1}{\sqrt{2\pi}} \int_{s}^\infty e^{-x^2/2}\,dx \sim
\exp(-s^2/2)/(\sqrt{2\pi} s).$$
For large $s$ (i.e., small $|S|$) this is dominated by $\exp(-s^2/2)$.
A similar statement holds for $T$. This shows that
Theorem~\ref{thm:isop} is nearly tight.

\section{The best asymptotic success rate in the $k$-star}

In this section we consider the NICD problem on the star. Let
$\Star_k$ denote the star graph on $k+1$ vertices and let $S_k$
denote its $k$ leaf vertices.  We shall study the same problem
considered in~\cite{MO03}; i.e., determining $\calM(\Star_k, \rho,
S_k)$.  Note that it was shown in that paper that the best
protocol in this case is always simple (i.e.,
all players should use the same function).\\

\noindent The following theorem determines rather accurately the
asymptotics of $\calM(\Star_k, \rho, S_k)$:

\begin{theorem} \label{thm:exactrate}
Fix $\rho \in (0,1]$ and let $\nu = \nu(\rho) = \frac{1}{\rho^2} -
1$. Then for $k \to \infty$,
\[
\calM(\Star_{k}, \rho, S_{k}) = \tilde{\Theta}\left(k^{-\nu}\right),
\]
where  $\tilde{\Theta}(\cdot)$ denotes asymptotics to within a
subpolynomial ($k^{o(1)}$) factor. The lower bound is achieved
asymptotically by the majority function $\maj_n$ with $n$
sufficiently large.
\end{theorem}

Note that if the corruption probability is very small (i.e.,
$\rho$ is close to 1), we obtain that the success rate only drops
off as a very mild function of $k$.

\begin{proof}[ of upper bound]
We know that all optimal protocols are simple, so assume all
players use the same balanced function $f : \bn \to \bits$.
Let $F_{-1} = f^{-1}(-1)$ and $F_1 = f^{-1}(1)$ be the sets where
$f$ obtains the values $-1$ and $1$ respectively.
The
center of the star gets a uniformly random string $x$, and then
independent $\rho$-correlated copies are given to the $k$ leaf
players.  Let $y$ denote a typical such copy.  The probability
that all players output $-1$ is thus $\E_x[\P[f(y) = -1 | x]^k]$. We
will show that this probability is $\tilde{O}(k^{-\nu})$. This
complete the proof since we can replace $f$ by $-f$ and get the
same bound for the probability that all players output $1$.


Suppose $\E_x[\P[f(y) = -1 | x]^k] \geq 2\delta$ for some
$\delta$; we will show $\delta$ must be small.  Define
\[
S = \{x : \Pr[f(y)= -1 \mid x]^k \geq \delta\}.
\]
By Markov's inequality we must have $|S| \geq \delta 2^{n}$. Now
on one hand, by the definition of $S$,
\begin{equation}\label{eqn:a1}
\Pr[y \in F_1 \mid x \in S] \leq 1 - \delta^{1/k}.
\end{equation}
On the other hand, applying Corollary~\ref{cor:isop} with $T =
F_1$ and $\alpha \leq 1/\log_2(1/\delta) < 1/\log(1/\delta)$
(since $|F_1| = \half
2^n$), we get
\begin{equation}\label{eqn:a2}
\Pr[y \in F_1 \mid x \in S] \geq
\delta^{(\log^{-1/2}(1/\delta) +
\rho)^2/(1-\rho^2)}.
\end{equation}
Combining~(\ref{eqn:a1}) and~(\ref{eqn:a2}) yields the desired
upper bound on $\delta$ in terms of $k$, $\delta \leq k^{-\nu +
o(1)}$ by the following calculations. We have
\[
1 - \delta^{1/k} \geq \delta^{(\log^{-1/2}(1/\delta) +
\rho)^2/(1-\rho^2)}.
\]
We want to show that the above inequality cannot hold if
\begin{equation} \label{eq:counter_delta}
\delta \geq \left(\frac{e^{c \sqrt{\log k}}}{k} \right)^{\nu},
\end{equation}
where $c = c(\rho)$ is some constant.
We will show that if $\delta$ satisfies (\ref{eq:counter_delta}) and
$c$ is sufficiently large then for all large $k$
\[
\delta^{1/k} +
\delta^{(\log^{-1/2}(1/\delta) + \rho)^2/(1-\rho^2)} > 1.
\]
Note first that
\begin{equation} \label{eq:left_delta}
\delta^{1/k} > \left( \frac{1}{k} \right)^{\frac{\nu}{k}} =
\exp \left(-\frac{\nu \log k}{k} \right) > 1 - \frac{\nu \log k}{k}.
\end{equation}
On the other hand,
\begin{equation} \label{eq:long_delta}
\delta^{(\log^{-1/2}(1/\delta) + \rho)^2/(1-\rho^2)} =
\delta^{- \log^{-1} \delta / (1-\rho^2) } ~\cdot~
\delta^{2 \rho  \log^{-1/2}(1/\delta) / (1-\rho^2)} ~\cdot~
  \delta^{\rho^2/(1-\rho^2)}.
\end{equation}
Note that
\[
  \delta^{\rho^2/(1-\rho^2)} = \delta^{1/\nu} \geq \frac{e^{c
  \sqrt{\log k}}}{k}
\]
and
\[
\delta^{2 \rho  \log^{-1/2}(1/\delta) / (1-\rho^2)} =
\exp\left( -\frac{2 \rho}{1 - \rho^2} \sqrt{\log(1/\delta)} \right)
\geq \exp\left( -\frac{2 \rho}{1 - \rho^2} \sqrt{\nu \log k} \right) .
\]
Finally,
\[
\delta^{- \log^{-1} \delta / (1-\rho^2) } =
\exp\left(-\frac{1}{1 - \rho^2} \right).
\]
Thus if $c = c(\rho)$ is sufficiently large then the left hand side of
(\ref{eq:long_delta}) is at least $\frac{\nu \log k}{k}$. This implies the
desired contradiction by (\ref{eq:counter_delta}) and (\ref{eq:left_delta}).
\end{proof}


\begin{proof}[ of lower bound]
We will analyze the protocol where all players use $\maj_n$, similarly
to the analysis of \cite{MO03}. Our analysis here is more careful resulting
in a tighter bound.

We begin by showing that the probability with which all
players agree if they use $\maj_n$, in the case of fixed $k$ and
$n \to \infty$, is:
\begin{equation}
\lim_{\substack{n \to \infty \\ n\text{ odd}}} \calP(\Star_{k},
\rho, n, S_{k}, \maj_n) = 2\nu^{1/2} (2\pi)^{(\nu
- 1)/2} \int_0^1 t^k I(t)^{\nu-1}\,dt, \label{eqn:theintegral}
\end{equation}
where $I = \phi \circ \Phi^{-1}$ is the so-called Gaussian
isoperimetric function, with $\phi(x) = (2 \pi)^{-1/2} \exp(-x^2/2)$
and $\Phi(x) = \int_{-\infty}^x \phi(t) dt$ the density and
distribution functions of a standard normal random variable respectively.

Apply Proposition~\ref{prop:normallimit}, with $X \sim N(0,1)$
representing $n^{-1/2}$ times the sum of the bits in the string at
the star's center, and $Y | X \sim N(\rho X, 1-\rho^2)$
representing $n^{-1/2}$ times the sum of the bits in a typical
leaf player's string. Thus as $n \to \infty$, the probability that
all players output $+1$ when using $\maj_n$ is precisely
\begin{equation*}
\infint \Phi\left(\frac{\rho\;x}{\sqrt{1-\rho^2}}\right)^k
\phi(x)\,dx = \infint \Phi\left(\nu^{-1/2} x \right)^k
\phi(x)\,dx.
\end{equation*}
Since $\maj_n$ is antisymmetric, the probability that all players
agree on $+1$ is the same as the probability they all agree on
$-1$.  Making the change of variables $t = \Phi(\nu^{-1/2} x)$, $x
= \nu^{1/2}\Phi^{-1}(t)$, $dx = \nu^{1/2} I(t)^{-1}\,dt$, we get
\begin{eqnarray*}
\lim_{\substack{n \to \infty \\ n\text{ odd}}} \calP(\Star_{k},
\rho, n, S_{k}, \maj_n) & = & 2\nu^{1/2}\int_0^1
\frac{t^k \phi(\nu^{1/2} \Phi^{-1}(t))}{ I(t) }\,dt \\
&=& 2\nu^{1/2} (2\pi)^{(\nu - 1)/2} \int_0^1 t^k
I(t)^{\nu-1}\,dt,
\end{eqnarray*}
as claimed.

We now estimate the integral in~(\ref{eqn:theintegral}).  It is
known (see, e.g.,~\cite{BG99}) that $I(s) \geq J(s(1-s))$, where
$J(s) = s\sqrt{\ln(1/s)}$.  We will forego the marginal
improvements given by taking the logarithmic term and simply use
the estimate $I(t) \geq t(1-t)$.  We then get
\begin{eqnarray*}
\int_0^1 t^k I(t)^{\nu-1}\,dt & \geq & \int_0^1 t^k
(t(1-t))^{\nu-1} \, dt \\
& = & \frac{\Gamma(\nu) \Gamma(k + \nu)}{\Gamma(k+2\nu)}
\qquad\quad\!\text{(\cite[6.2.1, 6.2.2]{AS72})} \\
& \geq & \Gamma(\nu) (k+2\nu)^{-\nu} \qquad \text{(Stirling
approximation).}
\end{eqnarray*}

Substituting this estimate into~(\ref{eqn:theintegral}) we get
$\lim_{n \to \infty} \calP(\Star_{k}, \rho, n, S_{k}, \maj_n) \geq
c(\nu) k^{-\nu}$ where $c(\nu) > 0$ depends only on $\rho$, as
desired.
\end{proof}

We remark that in the upper bound above
we have in effect proved the following theorem
regarding high norms of the Bonami-Beckner operator applied to
Boolean functions:
\begin{theorem} \label{thm:power} Let $f \colon \bn \to \{0,1\}$ and suppose $\E[f]
\leq 1/2$.  Then for any fixed $\rho \in (0,1]$, as $k \to
\infty$, $\|T_{\rho} f\|_k^k \leq k^{-\nu+ o(1)}$, where $\nu =
\frac{1}{\rho^2} - 1$.
\end{theorem}
Since we are trying to bound a high norm of $T_\rho f$ knowing the
norms of $f$, it would seem as though the usual Bonami-Beckner
inequality would be effective.  However this seems not to be the
case: a straightforward application yields
\begin{eqnarray*}
\|T_\rho f\|_k &\leq &\|f \|_{\rho^2(k-1)+1} = \E[f]^{1/(\rho^2(k-1)+1)}\\
\Rightarrow \qquad \|T_\rho f\|_k^k & \leq &
(1/2)^{k/(\rho^2(k-1)+1)} \approx (1/2)^{1/\rho^2},
\end{eqnarray*}
only a constant upper bound.

\section{The optimal protocol on the path}

In this section we prove the following theorem which gives a complete solution to the NICD
problem on a path. In this case, simple dictator protocols are the unique
optimal protocols, and any other simple protocol is exponentially worse
as a function of the number of players.

\begin{theorem}\label{thm:const_prot}
\begin{itemize}
\item
Let $\Path_k=\{v_0, v_1, \dots, v_k\} $ be the path graph of
length $k$, and let $S$ be any subset of $\Path_k$ of size at
least two. Then simple dictator protocols are the unique
optimal protocols for $\calP(\Path_k, \rho, n, S, (f_v))$.
In particular, if $S = \{v_{i_0}, \dots, v_{i_\ell}\}$ where $i_0
< i_1 < \cdots < i_\ell$, then we have
\[
\calM(\Path_k, \rho, S) = \prod_{j=1}^\ell \left(\frac12 + \frac12
\rho^{i_j - i_{j-1}}\right).
\]
\item
Moreover, for
every $\rho$ and $n$ there exists $c=c(\rho, n)< 1$ such that
if $S = \Path_k$ then for
any simple protocol $f$ which is not a dictator,
\[
\calP(\Path_k, \rho, n, S, f) \leq \calP(\Path_k, \rho, n, S, \Dicttr)
c^{|S|-1}
\]
where $\Dicttr$ denotes the dictator function.
\end{itemize}
\end{theorem}

\subsection{A bound on inhomogeneous Markov chains}

A crucial component of the proof of Theorem
\ref{thm:const_prot} is a bound on the probability that a Markov
chain stays inside certain sets. In this subsection, we derive
such a bound in a fairly general setting.
Moreover, we exactly characterize the cases in which the
bound is tight. This is a generalization of Theorem 9.2.7 in
\cite{AlonSpencer} and of results in \cite{AjtaiKS, AlonFWZ}.

Let us first recall some basic facts concerning reversible Markov
chains.  Consider an irreducible Markov chain on a finite set $S$.
We denote by $M=\big(m(x,y)\big)_{x,y \in S}$ the matrix of
transition probabilities of this chain, where $m(x,y)$ is a
probability to move in one step from $x$ to $y$. We will always assume
that $M$ is ergodic (i.e., irreducible and aperiodic).

The rule of the
chain can be expressed by the simple equation $\mu_1=\mu_0 M$,
where $\mu_0$ is a starting distribution on $S$ and
$\mu_1$ is a distribution obtained after one step of the Markov
chain (we think of both as row vectors). By definition, $\sum_y m(x,y)=1$. Therefore,
the largest eigenvalue of $M$ is $1$ and a corresponding right eigenvector has
all its coordinates equal to $1$.
Since $M$ is ergodic it has a unique (left/right) eigenvector corresponding
to an eigenvalue with absolute value $1$.
We denote the unique right eigenvector $(1,\ldots,1)^t$ by $\bf 1$.
We denote by $\pi$ the unique left eigenvector corresponding to the eigenvalue
$1$ whose coordinate sum is $1$. $\pi$
is the stationary distribution of the Markov chain. Since
we are dealing with a Markov chain whose distribution $\pi$ is not
necessarily uniform it will be convenient to work in $L^2(S,\pi)$.
In other words, for any two functions $f$ and $g$ on $S$ we define
the inner product $\ip{f, g} = \sum_{x\in S} \pi(x) f(x)g(x)$. The norm
of $f$ equals $\|f\|_2 = \sqrt{\ip{f,f}} = \sqrt{\sum_{x\in S} \pi(x) f^2(x)}$.

\begin{definition}
A transition matrix $M=\big(m(x,y)\big)_{x,y\in S}\,$ for a Markov chain is
reversible with respect to a probability distribution $\pi$ on $S$ if
$\pi(x)m(x,y)=\pi(y)m(y,x)$ holds for all $x,y$ in $S$.
\end{definition}

It is known that if $M$ is reversible with respect to $\pi$ then
$\pi$ is the stationary distribution of $M$. Moreover, the
corresponding operator taking $L^2(S,\pi)$ to itself defined by
$Mf(x)=\sum_y m(x,y) f(y)$ is self-adjoint, i.e.,
$\ip{Mf,g}=\ip{f,Mg}$ for all $f, g$. Thus, it follows that $M$
has a set of complete orthonormal (with respect to the inner
product defined above) eigenvectors with real eigenvalues.

\begin{definition}
If $M$ is reversible with respect to $\pi$ and $\lambda_1  \leq \ldots  \leq \lambda_{r-1} \leq \lambda_r=1$ are the
eigenvalues of $M$, then the {\em spectral
gap} of $M$ is defined to be $\delta=\min\big\{|-1-\lambda_1|,|1-\lambda_{r-1}|\big\}.$
\end{definition}

For transition matrices $M_1,M_2,\ldots$ on the same space $S$, we
can consider the time-inhomogeneous Markov chain which at time $0$
starts in some state (perhaps randomly) and then jumps using the
matrices $M_1,M_2,\ldots$ in this order. In this way,
$M_i$ will govern the jump from time $i-1$ to time $i$. We
write $I_{A}$ for the ($0$-$1$) indicator function of the set $A$
and $\pi_A$ for the function defined by $\pi_A(x) = I_A(x) \pi(x)$
for all $x$. Similarly, we define ${\pi}(A) = \sum_{x\in A}
\pi(x)$. The following theorem provides a tight estimate on the
probability that the inhomogeneous Markov chain stays inside
certain sets at every step.

\begin{theorem}\label{thm:aks}
Let $M_1,M_2,\ldots, M_k$ be ergodic transition
matrices on the state space $S$, all of which are reversible with
respect to the same probability measure $\pi$. Let $\delta_i>0$ be
the spectral gap of matrix $M_i$ and  let $A_0,A_1,\ldots,A_k$ be
subsets of $S$.
\begin{itemize}
\item
If $\{X_i\}_{i=0}^k$ denotes the
time-inhomogeneous Markov chain using the matrices
$M_1,M_2,\ldots, M_k$ and starting according to
distribution $\pi$, then
\begin{equation} \label{eq:aks_bd}
\Pr[X_i \in A_i\;\;\forall i=0 \ldots k] \le \sqrt{\pi(A_0)}
\sqrt{\pi(A_k)} \prod_{i=1}^{k} \Big[1 - \delta_{i} \left(
1 - \sqrt{\pi(A_{i-1})} \sqrt{\pi(A_{i})} \right) \Big].
\end{equation}
\item
Suppose we further assume that for all $i$, $\delta_i< 1$ and
that $\lambda^i_1 > -1+\delta_i$
($\lambda^i_1$ is the smallest eigenvalue for the $i$th chain).
Then equality in (\ref{eq:aks_bd})
holds if and only if the sets $A_i$ are the same set
$A$ and for all $i$ the function $I_A - \pi(A){\bf 1}$ is an
eigenfunction of $M_i$ corresponding to the eigenvalue
$1-\delta_i$.
\item
Finally, suppose even further that all the chains $M_i$ are
identical and that there is some set $A'$ such that equality holds
as above.  Then there exists a constant $c = c(M) < 1$ such that
for all sets $A$ for which a strict inequality holds, we have the
stronger inequality
\[
\Pr[X_i \in A_i\;\;\forall i=0 \ldots k] \le c^k \pi(A)
\prod_{i=1}^{k}\big[1 - \delta(1 - \pi(A)) \big].
\]
\end{itemize}
\end{theorem}

\begin{remark}
Notice that if all the sets $A_i$ have $\pi$-measure at most
$\sigma < 1$ and all the $M_i$'s have spectral gap at least
$\delta$, then the upper bound is bounded above by
$$
\sigma [\sigma+ (1-\delta)(1-\sigma)]^k.
$$
Hence, the above theorem generalizes the Theorem 9.2.7
in~\cite{AlonSpencer} and strengthens the estimate from \cite{AlonFWZ}.
\end{remark}

\subsection{Proof of Theorem \ref{thm:const_prot}}

If we look at the NICD process restricted to positions
$x_{i_0},x_{i_1},\ldots,x_{i_\ell}$, we obtain a
time-inhomogeneous Markov chain $\{X_j\}_{j=0}^\ell$ where $X_0$
is uniform on $\bn$ and the $\ell$ transition operators are powers
of the Bonami-Beckner operator, $T_\rho^{i_1-i_0},
T_\rho^{i_2-i_1}, \cdots, T_\rho^{i_\ell-i_{\ell-1}}$.
Equivalently, these operators are $T_{\rho^{i_1-i_0}}$,
$T_{\rho^{i_2-i_1}}$, \dots, $T_{\rho^{i_\ell-i_{\ell-1}}}$. It is
easy to see that the eigenvalues of $T_{\rho}$ are $1>\rho>
\rho^2> \dots>\rho^n$ and therefore its spectral gap is
$1-\rho$.  Now a protocol for the $\ell + 1$ players consists
simply of $\ell + 1$ subsets $A_0,\dots,A_\ell$ of $\bn$, where
$A_j$ is a set of strings in $\bn$ on which the $j$th player
outputs the bit $1$. Thus, each $A_j$ has size $2^{n-1}$, and the
success probability of this protocol is simply
\[
\Pr[X_i \in A_i\;\;\forall i=0 \ldots \ell] + \Pr[X_i \in
\bar{A}_i\;\;\forall i=0 \ldots \ell].
\]
But by Theorem~\ref{thm:aks} each summand is bounded by
$$
\frac{1}{2}\prod_{j=1}^\ell \left(\frac{1}{2} + \frac{\rho^{i_j -
i_{j-1}}}{2}\right),
$$
yielding our desired upper bound. It is easy to check that
this is precisely the success probability of a simple dictator protocol.

To complete the proof of the first part
it remains to show that every other protocol
does strictly worse.  By the second statement of
Theorem~\ref{thm:aks} (and the fact that the simple dictator
protocol achieves the upper bound in Theorem~\ref{thm:aks}), we
can first conclude that any optimal protocol is a simple protocol,
i.e., all the sets $A_j$ are identical. Let $A$ be the set
corresponding to any potentially optimal simple protocol. By
Theorem~\ref{thm:aks} again the function $I_A - (|A| 2^{-n}){\bf
1}=I_A -\frac{1}{2}{\bf 1}$ must be an eigenfunction of
$T_{\rho^r}$ for some $r$ corresponding to its second largest
eigenvalue $\rho^r$.  This implies that $f=2I_A -{\bf 1}$ must be
a balanced \emph{linear} function, $f(x) = \sum_{|S| = 1}
\hat{f}(S) x_S$. It is well known (see, e.g.,~\cite{FKN02}) that
the only such Boolean functions are dictators.
This completes the proof of the first part. The second part of the
theorem follows
immediately from the third part of Theorem \ref{thm:aks}
\qedsymb

\subsection{Inhomogeneous Markov chains}
In order to prove Theorem \ref{thm:aks} we need a lemma that provides
a bound on one step of the Markov chain.

\begin{lemma}\label{lem:inside_aks}
Let $M$ be an ergodic transition matrix for a
Markov chain on the set $S$ which is reversible with respect to
the probability measure $\pi$ and which has spectral gap $\delta>0$.
Let $A_1$ and $A_2$ be two subsets of $S$ and let $P_1$ and $P_2$
be the corresponding projection operators on $L^2(S,\pi)$ (i.e.,
$P_i f(x)=f(x)I_{A_i}(x)$ for every function $f$ on $S$).
Then
\[
\| P_2 M P_1\| \le 1 - \delta \left(1 - \sqrt{\pi(A_1)}
\sqrt{\pi(A_2)} \right),
\]
where the norm on the left is the operator norm for operators from
$L^2(S, \pi)$ into itself.

Further, suppose we assume that $\delta < 1$ and that $\lambda_1 >
-1+\delta$. Then equality holds above if and only if $A_1=A_2$ and
the function $I_{A_1} - \pi(A_1){\bf 1}$ is an eigenfunction of $M$
corresponding to $1-\delta$.
\end{lemma}

\begin{proof}
Let $e_1,\ldots,e_{r-1}, e_r={\bf 1}$ be an orthonormal basis of
right eigenvectors of $M$ with corresponding eigenvalues
$\lambda_1 \leq \ldots \leq \lambda_{r-1} \leq \lambda_r=1$.
For a function $f$ on $S$, denote by $\supp(f)=\{x \in S~|~ f(x)\not =
0\}$.
It is easy to see that
\[
\| P_2 M P_1\| =\sup\big\{|\ip{f_1,Mf_2}|:\|f_1\|_2=1,\|f_2\|_2=1,
\supp(f_1)\subseteq A_1, \supp(f_2)\subseteq A_2\big\}.
\]
Given such $f_1$ and $f_2$, expand them as
\[
f_1=\sum_{i=1}^r u_ie_i, \quad f_2=\sum_{i=1}^r v_ie_i
\]
and observe that for $j = 1, 2$,
\begin{equation} \label{eq:another_cs}
|\ip{f_j,{\bf 1}}|=|\ip{f_j,I_{A_j}}| \le \|f_j\|_2 \|I_{A_j}\|_2
=\sqrt{\pi(A_j)}.
\end{equation}
But now by the orthonormality of the $e_i$'s we have
\begin{eqnarray}\label{eq:inner_product}
|\ip{f_1,Mf_2}|&=& \left|\sum_{i=1}^r \lambda_i u_i v_i\right| \leq
\sum_{i=1}^r
|\lambda_i
u_i v_i| \leq | \ip{f_1,{\bf 1}} \ip{f_2,{\bf 1}}| + (1 - \delta) \sum_{i
\leq r-1}
|u_i v_i|\\
&\leq& \sqrt{\pi(A_1)} \sqrt{\pi(A_2)}+(1-\delta)\left(1-\sqrt{\pi(A_1)}
\sqrt{\pi(A_2)}\right)=1-\delta\left(1-\sqrt{\pi(A_1)}
\sqrt{\pi(A_2)}\right).\nonumber
\end{eqnarray}
Here we used that $\sum_i |u_i v_i| \le 1$ which follows from $f_1$ and
$f_2$ having
norm 1.

As for the second part of the lemma, if equality holds then all
the derived inequalities must be equalities. In particular, from the
inequality in (\ref{eq:another_cs}) it follows that for $j
= 1, 2$, $f_j = \pm \big(1/\sqrt{\pi(A_j)}\big) I_{A_j}$.
Since $\delta < 1$ is assumed it follows from the last
inequality in (\ref{eq:inner_product}) that
we must also have that $\sum_i |u_iv_i| = 1$ from
which we can conclude that $|u_i|=|v_i|$ for all $i$.
Since $-1 + \delta$ is not
an eigenvalue, for the last equality in \eqref{eq:inner_product}
to hold we must have that the only nonzero $u_i$'s (or $v_i$'s)
correspond to the eigenvalues $1$ and $1-\delta$.  Next, for the
first inequality in~\eqref{eq:inner_product} to hold, we must
have that $u=(u_1, \dots, u_n)=\pm v=(v_1, \dots, v_n)$ since $\lambda_i$
can only be $1$
or $1-\delta$
and $|u_i|=|v_i|$. This gives that $f_1=\pm f_2$ and therefore $A_1 =
A_2$.

Finally, we also get that $f_1 - \ip{f_1,{\bf 1}}{\bf 1}$ is an
eigenfunction of $M$ corresponding to the eigenvalue $1-\delta$.
To conclude the proof, note that if $A_1 = A_2$ and
$I_{A_1} - \pi(A_1){\bf 1}$ is an eigenfunction of $M$ corresponding to
$1-\delta$, then it is easy to see that when we take
$f_1=f_2=I_{A_1} - \pi(A_1){\bf 1}$, all inequalities in our proof become
equalities.
\end{proof}

\begin{proof}[ of Theorem \ref{thm:aks}]
Let $P_i$ denote the projection onto $A_i$, as in
Lemma~\ref{lem:inside_aks}. It is easy to see that
\[
\Pr[X_i \in A_i\;\;\forall i=0 \ldots k] =
\pi_{A_0}P_0 M_1 P_1 M_2 \cdots P_{k-1} M_k P_k I_{A_k}.
\]
Rewriting in terms of the inner product, this is equal to
\[
\ip{I_{A_0}, (P_0 M_1 P_1 M_2 \cdots P_{k-1} M_k P_k) I_{A_k}}.
\]
By Cauchy-Schwarz it is at most
\[
\|I_{A_0}\|_2 \|I_{A_k}\|_2 \| P_0 M_1 P_1 M_2 \cdots P_{k-1} M_k P_k\|,
\]
where the third factor is the norm of $P_0 M_1 P_1 M_2 \cdots P_{k-1} M_k
P_k$
as an operator from $L^2(S,\pi)$ to itself.
Since $P_i^2 = P_i$ (being a projection), this in turn is equal to
\[
\sqrt{\pi(A_0)} \sqrt{\pi(A_k)} \| (P_0 M_1 P_1) (P_1 M_2 P_2)\cdots
(P_{k-1} M_k P_k) \|.
\]
By Lemma~\ref{lem:inside_aks} we have that for all $i = 1, \dots,
k$
\[
\|P_{i-1} M_\ell P_{i}\| \leq 1 - \delta_i \Bigl(1 -
\sqrt{\pi(A_{i-1})}\sqrt{\pi(A_i)} \Bigr).
\]
Hence
\[
\Big\| \prod_{i = 1}^k (P_{i-1} M P_{i}) \Big\| \leq
\prod_{i=1}^{k} \left[1 - \delta_{i} \left( 1 -
\sqrt{\pi(A_{i-1})} \sqrt{\pi(A_{i})} \right) \right],
\]
and the first part of the theorem is complete.

For the second statement note that if we have equality, then we
must also have equality for each of the norms $\| P_{i-1} M_i P_{i}
\|$. This implies by Lemma~\ref{lem:inside_aks} that all the sets
$A_i$ are the same and that $I_{A_i} - \pi(A_i){\bf 1}$ is in the
$1-\delta_i$ eigenspace of $M_i$ for all $i$.   For the converse,
suppose on the other hand that $A_i = A$ for all $i$ and $I_A -
\pi(A){\bf 1}$ is in the $1-\delta_i$ eigenspace of $M_i$. Note that
\begin{eqnarray*}
P_{i-1} M_i P_i I_A &=& P_{i-1} M_i I_A = P_{i-1} M_i \big(\pi(A){\bf 1} +
(I_A - \pi(A){\bf 1})\big) =
P_{i-1} \big(\pi(A){\bf 1} + (1-\delta_i)(I_A -\pi(A){\bf 1})\big)\\
& =& \pi(A)I_A + (1 - \delta_i) I_A - (1 - \delta_i) \pi(A) I_A = \big(1 -
\delta_i(1 - \pi(A)\big)I_A.
\end{eqnarray*}
Since $P_i^2 = P_i$, we can use induction to show that
\[
\pi_{A_0}P_0 M_1 P_1 M_2 \cdots P_{k-1} M_k P_k I_{A_k}=
\pi_{A} \prod_{i=1}^k(P_{i-1} M_i P_i)\,I_{A} = \pi(A)
\prod_{i=1}^k \big(1 - \delta_{i} (1 - \pi(A)\big),
\]
completing the proof of the second statement.

In order to prove the third statement, it suffices to note that if
$A$ does not achieve equality and $P$ is the corresponding
projection, then $\|PMP\| < 1 - \delta(1 - \pi(A))$.
\end{proof}

\section{NICD on general trees}
In this section we give some results for the NICD problem on
general trees.
First we observe that the following statement follows easily from
the proof of Theorem 1.3 in \cite{MO03}:

\begin{theorem}
\label{small-S}
For any NICD instance $(T, \rho, n, S)$ in which $|S| = 2$ or $|S|
= 3$ the simple dictator protocols constitute all optimal
protocols.
\end{theorem}

\subsection{Example with no simple optimal protocols}
It appears that the problem of NICD in general
is quite difficult. In particular, using
Theorem~\ref{thm:const_prot} we show that there are instances for
which there is no simple optimal protocol. Note the contrast with the
case of stars where it is proven in \cite{MO03} that there is always a
simple optimal protocol.

\begin{theorem}\label{thm:no_const_prot}
There exists an instance $(T,\rho, n, S)$ for which there is no
simple optimal protocol. In fact, given any $\rho$ and any $n \geq
4$, there are integers $k_1$ and $k_2$, such that if $T$ is a
$k_1$-leaf star together with a path of length $k_2$ coming out of
the center of the star (see Figure~\ref{fig:starpath})
and $S$ is the full vertex set of $T$,
then this instance has no simple optimal protocol.
\end{theorem}
\begin{proof}
Fix $\rho$ and $n \geq 4$.
Recall that we write $\eps = \half - \half \rho$ and let
$\mathrm{Bin}(3, \eps)$ be a binomially distributed random
variable with parameters $3$ and $\eps$. As was observed
in~\cite{MO03},
\[
\calP(\Star_{k}, \rho, n, S_{k}, \maj_3) \geq \frac18
\Pr[\mathrm{Bin}(3, \eps) \leq 1]^k.
\]
To see this, note that with probability $1/8$ the center of the
star gets the string $(1, 1, 1)$.  Since \linebreak
$\Pr[\mathrm{Bin}(3, \eps) \leq 1]= (1 - \eps)^2(1+2\eps)
> 1 - \eps$ for all $\epsilon<1/2$, we can pick $k_1$ large enough
so that
$$\calP(\Star_{k_1}, \rho, n, S_{k_1}, \maj_3) \geq 8(1 -
\eps)^{k_1}.$$

Next, by the last statement in Theorem~\ref{thm:aks},
there exists $c_2=c_2(\rho, n) > 1$ such that for all balanced
non-dictator functions $f$ on $n$ bits
$$
{\cal P}(\Path_k,\rho,n,\Path_k,\Dicttr) \ge {\cal P}(\Path_k, \rho, n, \Path_k, f)
c_2^k.
$$
Choose $k_2$ large enough so that
$$
(1 - \eps)^{k_1} c_2^{k_2} > 1.
$$

\begin{figure}[h]
\center{
 \epsfxsize=2in\epsfbox{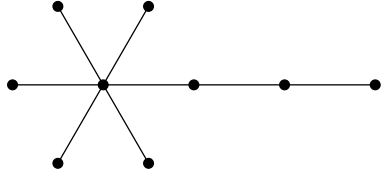}}
 \caption{The graph $T$ with $k_1=5$ and $k_2=3$}
 \label{fig:starpath}
\end{figure}

Now let $T$ be the graph consisting of a star with $k_1$ leaves
and a path of length $k_2$ coming out of its center (see
Figure~\ref{fig:starpath}), and let $S=V(T)$.  We claim that the
NICD instance $(T,\rho,n,S)$ has no simple optimal protocol.  We
first observe that if it did, this protocol would have to be
$\Dicttr$, i.e., ${\cal P}(T,\rho,n,S,f) < {\cal
P}(T,\rho,n,S,\Dicttr)$ for all simple protocols $f$ which are not
equivalent to dictator.  This is because the quantity on the right
is $(1 - \eps)^{k_1+k_2}$ and the quantity on the left is at most
${\cal P}(\Path_{k_2},\rho,n,\Path_{k_2},f)$ which in turn by definition
of $c_2$ is at most $(1 - \eps)^{k_2} / c_2^{k_2}$. This is
strictly less than $(1 - \eps)^{k_1+k_2}$ by the choice of $k_2$.

To complete the proof it remains to show that $\Dicttr$ is not an
optimal protocol.  Consider the protocol where $k_2$ vertices on
the path (including the star's center) use the dictator
$\Dicttr$ on the first bit and the $k_1$ leaves of the star use the protocol $\maj_3$
on the last three out of $n$ bits.  Since $n \geq 4$, these
vertices use completely independent bits from those that vertices
on the path are using.  We will show that this protocol, which we
call $f$, does better than $\Dicttr$.

Let $A$ be the event that all vertices on the path have their
first bit being 1. Let $B$ be the event that each of the $k_1$
leaf vertices of the star have 1 as the majority of their last 3
bits.  Note that $P(A) = \half (1 - \eps)^{k_2}$ and that, by
definition of $k_1$,  $P(B) \ge 4 (1 - \eps)^{k_1}$. Now the
protocol $f$ succeeds if both $A$ and $B$ occur.  Since $A$ and
$B$ are independent (as distinct bits are used), $f$ succeeds with
probability at least $2 (1 - \eps)^{k_2} (1 - \eps)^{k_1}$ which
is twice the probability that the dictator protocol succeeds.
\end{proof}

\begin{remark}
It was not necessary to use the last 3 bits for the $k_1$
vertices; we could have used the first 3 (and had $n = 3$).  Then
$A$ and $B$ would not be independent but it is easy to show (using
the FKG inequality) that $A$ and $B$
would then be positively correlated which is all that is needed.
\end{remark}

\subsection{Optimal monotone protocols}
Next, we present some general statements about what optimal
protocols must look like.  Using discrete symmetrization together
with the FKG inequality we prove the following theorem, which
extends one of the results in~\cite{MO03} from the case of the star
to the case of general trees.

\begin{theorem}\label{thm:optimal_monotone}
For all NICD instances on trees, there is an optimal protocol in
which all players use a monotone function.
\end{theorem}

One of the tools that we need to prove Theorem
\ref{thm:optimal_monotone} is the correlation inequality obtained
by Fortuin et al.~\cite{FKG} which is usually called the FKG
inequality. We first recall some basic definitions.

Let $D$ be a finite linearly ordered set.
Given two strings $x, y$ in $D^m$ we write $x
\leq y$ iff $x_i\leq y_i$ for all indices $1\leq i \leq m$. We
denote by $x \vee y$ and $x \wedge y$ two strings whose $i$th
coordinates are $\max(x_i,y_i)$ and $\min(x_i,y_i)$ respectively.
A probability measure $\mu: D^m \rightarrow \Reals^{\geq 0}$ is
called {\em log-supermodular} if
\begin{equation} \label{log}
\mu(\eta) \mu(\delta) \leq \mu(\eta \vee \delta)
\mu(\eta \wedge \delta)
\end{equation}
for all $\eta, \delta \in D^m$. If $\mu$ satisfies
(\ref{log}) we will also say that $\mu$ satisfies the
FKG condition. A subset $A \subseteq D^m$ is {\em
increasing} if whenever $x \in A$ and $x \leq y$ then also $y \in
A$. Similarly, $A$ is {\em decreasing} if $x \in A$ and $y \leq x$
imply that $y \in A$. Finally, the measure of $A$ is
$\mu(A)=\sum_{x \in A}\mu(x)$. The following well known fact is a
special case of the FKG inequality.

\begin{proposition}
\label{FKG-ineq} Let $\mu: \bm \rightarrow \Reals^{\geq 0}$ be
a log-supermodular probability measure on the discrete cube. If $A$
and $B$ are two increasing subsets of $\bm$ and $C$ is a
decreasing subset then
$$\mu(A\cap B) \geq \mu(A) \cdot \mu(B)~~~~ \mbox{and} ~~~~
\mu(A\cap C) \leq \mu(A) \cdot \mu(C).$$
\end{proposition}

It is known that in order to
prove that $\mu$ satisfies the FKG lattice condition,
it suffices to check this for
``smallest boxes'' in the lattice, i.e., for $\eta$ and $\delta$ that agree
at all but two locations. For completeness we prove this here.
\begin{lemma} \label{lem:small_box}
Let $\mu$ be a measure with full support. Then $\mu$ satisfies the FKG
condition (\ref{log}) if and only if it satisfies (\ref{log})
for all $\eta$ and $\delta$ that agree at all but two locations.
\end{lemma}

\begin{proof}
We will prove the non-trivial direction
by induction on $d = d(\eta,\delta)$, the Hamming
distance between $\eta$ and $\delta$. The cases where $d(\eta,\delta)
\leq 2$ follow from the assumption. The proof will proceed by induction on
$d$. Let $d = d(\eta,\delta) \geq 3$ and assume the claim holds for
all smaller $d$. We can partition the set of coordinates into $3$
subsets $I_{=}, I_{\{\eta > \delta\}}$ and $I_{\{\eta < \delta\}}$, where
$\eta$ and $\delta$ agree, where $\eta > \delta$ and where $\eta < \delta$
respectively.
Without loss of generality $|I_{\{\eta > \delta\}}| \geq 2$.
Let $i \in I_{\{\eta > \delta\}}$ and let $\eta'$ be obtained from
$\eta$ by setting $\eta'_i = \delta_i$ and letting $\eta'_j = \eta_j$
otherwise. Then since $\eta' \wedge \delta = \eta \wedge \delta$,
\[
\frac{\mu(\eta \wedge \delta) \mu(\eta \vee
  \delta)}{\mu(\eta)\mu(\delta)} =
\left(\frac{\mu(\eta' \wedge \delta) \mu(\eta' \vee
  \delta)}{\mu(\eta')\mu(\delta)} \right) \times
\left( \frac{\mu(\eta') \mu(\eta \vee \delta)}{\mu(\eta)\mu(\eta' \vee
  \delta)} \right).
\]
The first factor is $\geq 1$ by the induction hypothesis since
$d(\delta,\eta') = d(\delta,\eta)-1$. Note that $\eta' = \eta \wedge
(\eta' \vee \delta)$ and $\eta \vee \delta = \eta \vee (\eta' \vee
\delta)$ and $d(\eta',\eta \vee \delta) = 1 + |I_{\{\eta < \delta\}}| <
d$. Therefore by induction, the second term is also $\geq 1$.
\end{proof}

The above tools together with symmetrization now allow us to prove
Theorem~\ref{thm:optimal_monotone}.

\begin{proof}[ of Theorem~\ref{thm:optimal_monotone}]
Recall that we have a tree $T$ with $m$ vertices, $0<\rho < 1$, and a probability
measure $\P$ on $\alpha \in \bits^{V(T)}$ which is defined by
\[
\P(\alpha) = \half (\half+\half\rho)^{A(\alpha)} (\half -
\half\rho)^{B(\alpha)},
\]
where $A(\alpha)$ is the number of pairs of neighbors where
$\alpha$ agrees and $B(\alpha)$ is the number of pairs of
neighbors where $\alpha$ disagrees. To use Proposition
\ref{FKG-ineq} we need to show that $\P$ is a log-supermodular
probability measure.

Note that (\ref{log}) holds  trivially if $\alpha \leq \beta$  or
$\beta \leq \alpha$. Thus it suffices to consider the
case where there are two vertices $u,v$ of $T$ on which
$\alpha$ and $\beta$ disagree and that
$\alpha_v=\beta_u=1$ and $\alpha_u=\beta_v=-1$. If these vertices
are not neighbors then by definition of $\P$ we have that
$\P(\alpha)\P(\beta)=\P(\alpha\vee \beta)\P(\alpha \wedge \beta)$.
Similarly, if $u$ is a neighbor of $v$ in $T$, then one can easily
check that
$$\frac{\P(\alpha)\P(\beta)}{\P(\alpha\vee \beta)\P(\alpha \wedge \beta)}=
\left(\frac{1-\rho}{1+\rho}\right)^2 \leq 1.$$ Hence we conclude
that measure $\P$ is log-supermodular.

Let $f_1,\dots,f_k$ be the
functions used by the parties at nodes $S = \{v_1,\dots,v_k\}$. We
will shift the functions in the sense of Kleitman's monotone
``down-shifting''~\cite{Kle66}.
Namely, define functions $g_1,\dots,g_k$
as follows: If $f_i(-1,x_2,\dots,x_n) = f_i(1,x_2,\dots,x_n)$ then
we set
$$
g_i(-1,x_2,\dots,x_n) = g_i(1,x_2,\dots,x_n) =
f_i(-1,x_2,\dots,x_n) = f_i(1,x_2,\dots,x_n).
$$
Otherwise, we set $g_i(-1,x_2,\dots,x_n) = -1$ and
$g_i(1,x_2,\dots,x_n) = 1$. We claim that the agreement
probability for the $g_i$'s is at least the agreement probability
for the $f_i$'s. Repeating this argument for all bit locations
will prove that there exists an optimal protocol for which all
functions are monotone.

To prove the claim we condition on the value of $x_2,\dots,x_n$ at
all the nodes $v_i$ and let $ \alpha_i$ be the remaining bit at
$v_i$. For simplicity we will denote the functions of this bit by
$f_i$ and $g_i$. Note that if there exist $i$ and $j$ such that
$f_i(-1) = f_i(1) = -1$ and $f_j(-1) = f_j(1) = 1$, then the
agreement probability for both $f$ and $g$ is $0$.

It therefore
remains to consider the case where there exists a subset $S'
\subset S$ such that $f_i(-1) = f_i(1) = 1$ for all $i \in S'$ and
$f_i(-1) \neq f_i(1)$ for all $i \in U=S \setminus S'$ (the case
where $f_i(-1) = f_i(1) = -1$ for all $i \in S'$ can be treated
similarly and the case where for all functions $f_i(-1) \neq
f_i(1)$ may be decomposed into the above two events where $S' =
\emptyset$). Note that in this case the agreement probability for
the $g$'s is nothing but $\P(\alpha_i = 1 : i \in U)$ while the
agreement probability for the $f$'s is $\P(\alpha_i =\tau_i : i
\in U)$, where $\tau_i = -1$ if $f_i(-1)=1$ and $\tau_i=1$
otherwise.

Let $U' \subseteq U$ be the set of indices $i$ such
that $\tau_i = -1$ and let $U^{''}=\{i \in U ~|~ \tau_i = 1\}$.
Let $A$ be the set of strings in $\bm$ with $\alpha_i=1$ for all $i
\in U'$, let $B$ be the set of strings with $\alpha_i=1$ for all $i
\in U^{''}$ and let $C$ be the set of strings with $\alpha_i=-1$ for
all $i \in U'$. Note that $A, B$ are increasing sets and $C$ is
decreasing. Also, since our distribution is symmetric, it is easy
to see that $\P(A)=\P(C)$. Therefore, by the FKG inequality, the
agreement probability for the $g$'s, namely
$$\P(A\cap B) \geq \P(A)\cdot\P(B)=\P(C)\cdot\P(B)\geq \P(C\cap B),$$  is at least as
large as for the $f$'s.
\end{proof}

\begin{remark} The last step in the proof above may be
  replaced by a more direct calculation showing that in
  fact we have strict inequality unless the sets $U',U''$ are empty.
  This is similar to the monotonicity proof in \cite{MO03}.
  This implies that every optimal protocol must
consist of monotone functions (in general, it may be monotone increasing in
some coordinates and monotone decreasing in the other coordinates).
\end{remark}

\begin{remark} The above proof works in a much more general setup than just our
tree-indexed Markov chain case. One can take any FKG measure on $\bits^m$ with all
marginals having mean $0$, take $n$ independent copies of this and define
everything analogously in this more general framework. The proof of
Theorem~\ref{thm:optimal_monotone} extends to this context.
\end{remark}

\subsection{Monotonicity in the number of parties}
Our last theorem yields a certain monotonicity when comparing the
simple dictator protocol $\Dicttr$ and the simple protocol
$\maj_r$, which is majority on the first $r$ bits.  The result is not very
strong -- it is interesting mainly because it allows to compare protocols
behavior for different number of parties. It shows that if $\maj_r$ is a
better protocol than dictatorship for $k_1$ parties on the star,
then it is also better than dictatorship for $k_2$ parties if $k_2 > k_1$.

\begin{theorem}\label{thm:monotonicity}
Fix $\rho$ and $n$ and suppose $k_1$ and $r$ are such that
\[
\calP(\Star_{k_1}, \rho, n, \Star_{k_1}, \maj_r) \geq
(>)\;\calP(\Star_{k_1}, \rho, n, \Star_{k_1}, \Dicttr).
\]
Then for all $k_2 > k_1$,
\[
\calP(\Star_{k_2}, \rho, n, \Star_{k_2}, \maj_r) \geq
(>)\;\calP(\Star_{k_2}, \rho, n, \Star_{k_2}, \Dicttr).
\]
\end{theorem}

Note that it suffices to prove the theorem assuming $r=n$.
In order to prove the theorem,
we first recall the notion of stochastic domination. If
$\eta,\delta\in\{0,1,\ldots,n\}^I$, write $\eta\preceq\delta$ if
$\eta_i\le\delta_i$ for all $i\in I$. If  $\nu$ and $\mu$ are two
probability measures on $\{0,1,\ldots,n\}^I$, we say $\mu$ {\em
stochastically dominates}  $\nu$, written $\nu \preceq \mu$, if
there exists a probability measure $m$ on $\{0,1,\ldots,n\}^I
\times \{0,1,\ldots,n\}^I$ whose first and second marginals are
respectively $\nu$ and $\mu$ and such that $m$ is supported on
$\{(\eta,\delta): \eta\preceq\delta\}$.

Fix $\rho$, $n \geq 3$, and any tree $T$.  Let our tree-indexed
Markov chain be $\{x_v\}_{v\in T}$, where $x_v \in \bn$ for each
$v\in T$.  Let $A\subseteq \bn$ be the strings which have a
majority of 1's.  Let $X_v$ denote the number of 1's in $x_v$.

Given $S\subseteq T$, let $\mu_S$ be the conditional distribution
of $\{X_v\}_{v\in T}$ given $\cap_{v\in S}\{x_v \in A\}$
($=\cap_{v\in S}\{X_v \ge n/2\}$). The following lemma is key and
might be of interest in itself. It can be used to prove (perhaps
less natural) results analogous to Theorem~\ref{thm:monotonicity}
for general trees.  Its proof will be given later.

\begin{lemma} \label{lem:5.1}  In the above setup, if $S_1 \subseteq S_2\subseteq T$, we have
\[
\mu_{S_1}\preceq \mu_{S_2}.
\]
\end{lemma}

Before proving the lemma or showing how it implies
Theorem~\ref{thm:monotonicity}, a few remarks are in order.

\begin{itemize}
\item
Note that if $\{x_k\}$ is a Markov chain
on $\bn$ with transition matrix $T_\rho$, then if we let $X_k$ be
the number of 1's in $x_k$, then $\{X_k\}$ is also a Markov chain
on the state space $\{0,1,\dots,n\}$ (although it is certainly not true in
general that a function of a Markov chain is a Markov chain.)  In
this way, with a slight abuse of notation, we can think of
$T_\rho$ as a transition matrix for $\{X_k\}$ as well as for
$\{x_k\}$.  In particular, given a probability distribution $\mu$
on $\{0, 1, \dots, n\}$ we will write $\mu T_\rho$ for the probability
measure on $\{0,1,\dots,n\}$ given by one step of the Markov
chain.
\item
We next recall the easy fact that the Markov chain $T_\rho$ on
$\bn$ is \emph{attractive} meaning that if $\nu$ and $\mu$ are
probability measures on $\bn$ with $\nu \preceq \mu$, then it
follows that $\nu T_\rho  \preceq \mu T_\rho $.
The same is true for the Markov chain
$\{X_k\}$ on $\{0,1,\dots,n\}$.
\end{itemize}
Along with these observations, Lemma~\ref{lem:5.1} is enough to
prove Theorem~\ref{thm:monotonicity}:

\begin{proof}
Let $v_0,v_1,\dots,v_k$ be the vertices of $\Star_k$, where $v_0$ is the center.
Clearly,
${\cal P}(\Star_k, \rho, \Star_k,\Dicttr)=(\half + \half \rho)^k$.
On the other hand, a little thought reveals that
\[
{\cal P}(\Star_k,\rho, n, \Star_k,\maj_n)=
\prod_{\ell=0}^{k-1}(\mu_{v_0,\ldots,v_\ell}\mid_{v_0} T_\rho)(A),
\]
where $\nu\mid_v$ means the distribution of $\nu$ restricted to
the location $v$ (recall that $A\subseteq \bn$ is the strings which have a
majority of 1's).
By Lemma~\ref{lem:5.1} and the attractivity of
the process, the terms $(\mu_{v_0,\ldots,v_\ell}\mid_{v_0}
T_\rho)(A)$ (which do not depend on $k$ as long as $\ell \le k$)
are nondecreasing in $\ell$. Therefore if
\[
{\cal P}(\Star_k,\rho, n, \Star_k,\maj_n) \geq (>) (\half + \half \rho)^k,
\]
then
$(\mu_{v_0,\ldots,v_{k-1}}\mid_{v_0} T_\rho)(A) \geq (>) \half + \half \rho$
which implies in turn that for every $k' \geq k$,
$(\mu_{v_0,\ldots,v_{k'-1}}\mid_{v_0} T_\rho)(A) \geq (>) \half + \half \rho$
and thus for all $k' > k$
\[
{\cal P}(\Star_{k'},\rho, n, \Star_{k'},\maj_n) \geq (>) (\half + \half \rho)^{k'}.
\]

\end{proof}

\bigskip

Before proving Lemma~\ref{lem:5.1}, we recall the definition of
\emph{positive associativity}. If $\mu$ is a probability measure
on $\{0,1,\dots,n\}^I$, $\mu$ is said to be \emph{positively
associated} if any two functions on $\{0,1,\dots,n\}^I$ which are
increasing in each coordinate are positively correlated. This is
equivalent to the fact that if $B\subseteq \{0,1,\dots,n\}^I$ is
an upset, then $\mu$ conditioned on $B$ is stochastically larger
than $\mu$.

\begin{proof}[ of Lemma~\ref{lem:5.1}]
It suffices to prove this when $S_2$ is $S_1$ plus an extra
vertex $z$. We claim that for any set $S$, $\mu_S$ is positively
associated. Given this claim, we form $\mu_{S_2}$ by first
conditioning on $\cap_{v\in S_1}\{x_v\in A\}$, giving us the
measure $\mu_{S_1}$, and then further conditioning on $x_z\in A$.
By the claim, $\mu_{S_1}$ is positively associated and hence the
last further conditioning on $X_z\in A$ stochastically increases
the measure, giving $\mu_{S_1}\preceq \mu_{S_2}$.

To prove the claim that $\mu_S$ is positively associated, we first
claim that the distribution of $\{X_v\}_{v\in T}$, which is just a
probability measure on $\{0,1,\dots,n\}^T$, satisfies the FKG lattice condition (\ref{log}).

Assuming the FKG condition holds
for $\{X_v\}_{v \in T}$, it is easy to see that the
same inequality holds when we condition on the sublattice
$\cap_{v\in S}\{X_v\ge n/2\}$ (it is crucial here that the set
$\cap_{v\in S}\{X_v\ge n/2\}$ is a sublattice meaning that
$\eta,\delta$ being in this set implies that $\eta\,\,  \vee \,\,
\delta$ and $\eta\,\,  \wedge\,\, \delta$ are also in this set).

The FKG theorem, which says that the FKG lattice
condition (for any distributive lattice) implies positive
association, can now be applied to this conditioned measure to
conclude that the conditioned measure has positive association, as desired.

Finally, by Lemma \ref{lem:small_box}, in order to
prove that $P$ satisfies the FKG lattice condition,
it is enough to check this for
``smallest boxes'' in the lattice, i.e., for $\eta$ and $\delta$ that agree
at all but two locations.
If these two locations are not neighbors (i.e., two leaves), it is easy to
check that we have equality. If they are neighbors, it easily
comes down to checking that if $a > b$ and $c > d$, then
$$
P(X_1=c|X_0=a) P(X_1=d|X_0=b) \ge P(X_1=d|X_0=a) P(X_1=c|X_0=b)
$$
where $\{X_0,X_1\}$ is the distribution of our Markov chain on
$\{0,1,\dots,n\}$ restricted to two consecutive times. It is
straightforward to check that for $\rho \in (0,1)$, the above
Markov chain can be embedded into a continuous time Markov chain
on $\{0,1,\dots,n\}$ which only takes steps of size $1$. Hence the
last claim is a special case of Lemma \ref{lem:reflection}
\end{proof}

\begin{lemma} \label{lem:reflection}

If $\{X_t\}$ is a continuous time Markov chain on $\{0,1,\dots,n\}$
which only takes steps of size 1, then if $a > b$ and $c > d$, it follows that
$$
\P(X_1=c~|~X_0=a) ~\P(X_1=d~|~X_0=b) \ge \P(X_1=d~|~X_0=a) ~\P(X_1=c~|~X_0=b).
$$
(Of course, by time scaling, $X_1$ can be replaced by any time $X_t$.)
\end{lemma}

\begin{proof}
Let $R_{a,c}$ be the set of all possible realizations of our
Markov chain during $[0,1]$ starting from $a$ and ending in $c$.
Define $R_{a,d}$, $R_{b,c}$ and $R_{b,d}$ analogously. Letting
$P_x$ denote the measure on paths starting from $x$, we need to
show that
$$
P_a(R_{a,c})P_b(R_{b,d}) \ge P_a(R_{a,d}) P_b(R_{b,c})
$$
or equivalently that
$$
P_a\times P_b [R_{a,c}\times R_{b,d}] \ge P_a\times P_b [R_{a,d}\times R_{b,c}]
$$
We do this by giving a measure preserving injection from
$R_{a,d}\times R_{b,c}$ to $R_{a,c}\times R_{b,d}$. We can ignore
pairs of paths where there is a jump in both paths at the same time since these
have $P_a\times P_b$ measure 0. Given a pair of paths in
$R_{a,d}\times R_{b,c}$, we can switch the paths after their first
meeting time. It is clear that this gives an injection from
$R_{a,d}\times R_{b,c}$ to $R_{a,c}\times R_{b,d}$ and the Markov
property guarantees that this injection is measure preserving,
completing the proof.
\end{proof}


\section{Conclusions and open questions}

In this paper we have exactly analyzed the NICD problem on the path
and asymptotically analyzed the NICD problem on the star.  However,
we have seen that results on more complicated trees may be hard to
come by. Many problems are still open. We list a few:
\begin{itemize}
\item
Is it true that for every tree NICD instance,
     there is an optimal protocol in which each player uses some
     majority rule? This question was already raised in~\cite{MO03}
     for the special case of the star.
\item
Our analysis for the star is quite tight. However, one can ask for more.
In particular, what is the best bound that can be obtained on
\[
r_k = \frac{\calM(\Star_k, \rho, S_k)}
{\lim_{\substack{n \to \infty \\ n\text{ odd}}} \calP(\Star_{k},
\rho, n, S_{k}, \maj_n)}
\]
for fixed value of $\rho$.
Our results show that $r_k = k^{o(1)}$. Is it true that
$\lim_{k \to \infty} r_k = 1$?
\item
Finally, we would like to find more applications of the reverse
Bonami-Beckner inequality in computer science and combinatorics.
\end{itemize}

\section{Acknowledgments}

Thanks to David Aldous, Christer Borell, Svante Janson, Yuval Peres, and Oded
Schramm for helpful discussions.

\bibliographystyle{abbrv}
\bibliography{beckner}

\appendix

\section{Proof of the reverse Bonami-Beckner inequality}\label{sec:proof_rev_beckner}
Borell's proof of the reverse Bonami-Beckner
inequality~\cite{Bor82} follows the same lines as the traditional
proofs of the usual Bonami-Beckner inequality~\cite{Bon70,Bec75}.  Namely, he proves the result in the case $n
= 1$ (i.e., the ``two-point inequality'') and then shows that this
can be tensored to produce the full theorem.  The usual proof of
the tensoring is easily modified by replacing Minkowski's
inequality with the reverse Minkowski inequality~\cite[Theorem 24]{HLP34}.
Hence, it is enough to consider functions $f : \bits \to \Reals^{\geq 0}$
(i.e., $n=1$).
By monotonicity of norms, it
suffices to prove the inequality in the case that $\rho =
(1-p)^{1/2}/(1-q)^{1/2}$; i.e., $\rho^2 = (1-p)/(1-q)$.
Finally, it turns out that it suffices to
consider the case where $0 < q < p < 1$ (see Lemma \ref{lem:duality}).
\begin{lemma} Let $f : \bits \to \Reals^{\geq 0}$ be a nonnegative
function, $0 < q < p < 1$, and $\rho^2 = (1-p)/(1-q)$.  Then
$\|T_\rho f \|_q \geq \|f\|_p$.
\end{lemma}
\begin{proof}[ (Borell)]  If $f$ is identically zero the lemma is
trivial.  Otherwise, using homogeneity we may assume that $f(x) =
1 + ax$ for some $a \in [-1,1]$.  We shall consider only the case
$a \in (-1,1)$; the result at the endpoints follows by
continuity.  Note that $T_\rho f (x) = 1 + \rho a x$.

Using the Taylor series expansion for $(1+a)^q$ around 1, we get

\begin{eqnarray}
\|T_{\rho} f\|_q^q &=& \frac{1}{2}\left( (1 + a \rho)^q + ( 1 - a
\rho)^q \right) = \frac{1}{2} \left((1 + \sum_{n=1}^{\infty}
\binom{q}{n} a^n \rho^n) +
      (1 + \sum_{n=1}^{\infty} \binom{q}{n} (-a)^n \rho^n) \right)
\nonumber\\ &=& 1 + \sum_{n=1}^{\infty} \binom{q}{2n} a^{2n}
\rho^{2n}. \label{eq:Tq}
\end{eqnarray}
(Absolute convergence for $|a| < 1$ lets us rearrange the series.)
 Since $p > q$, it holds for all $x > -1$ that $(1 + x)^{p/q} \geq
1 + px/q$. In particular, from (\ref{eq:Tq}) we obtain that
\begin{equation} \label{eq:Tp}
\|T_{\rho} f\|_q^p = \left( 1 + \sum_{n=1}^{\infty} \binom{q}{2n}
a^{2n} \rho^{2n} \right)^{p/q} \geq 1 + \sum_{n=1}^{\infty}
\frac{p}{q} \binom{q}{2n} a^{2n} \rho^{2n}.
\end{equation}
Similarly to (\ref{eq:Tq}) we can write
\begin{equation} \label{eq:fp}
\|f\|_p^p =
1 + \sum_{n=1}^{\infty} \binom{p}{2n} a^{2n}.
\end{equation}
 From (\ref{eq:Tp}) and (\ref{eq:fp}) we see that in order to prove
the theorem it suffices to show that for all $n \geq 1$
\begin{equation} \label{eq:comp_coeff}
\frac{p}{q} \binom{q}{2n} \rho^{2n} \geq \binom{p}{2n}.
\end{equation}
Simplifying (\ref{eq:comp_coeff}) we see the inequality
\[
(q-1)\cdots (q-2n+1) \rho^{2n} \geq (p-1) \cdots (p-2n+1),
\]
which is equivalent in turn to
\begin{equation} \label{eq:comp_coeff2}
(1-q)\cdots (2n-1-q) \rho^{2n} \leq (1-p) \cdots (2n-1-p).
\end{equation}
Note that we have $(1-p) = (1-q) \rho^2$. Inequality
(\ref{eq:comp_coeff}) would follow if we could show that for all $m
\geq 2$ it holds that $\rho (m - q) \leq (m - p)$. Taking the square
and recalling that $\rho^2 = (1-p)/(1-q)$ we obtain the inequality
\[
(1-p) (m-q)^2 \leq (m-p)^2 (1-q),
\]
which is equivalent to
\[
m^2 - 2 m + p + q - pq \geq 0.
\]
The last inequality holds for all $m \geq 2$ thus completing the
proof.
\end{proof}

\bigskip

We also prove the two-function version promised in Section~\ref{sec:wherearewe}.
Recall first the reverse H\"{o}lder inequality \cite[Theorem 13]{HLP34}
for discrete measure spaces:
\begin{lemma}
Let $f$ and $g$ be nonnegative functions  and suppose $1/p + 1/p' = 1$,
where $p < 1$ ($p' = 0$ if $p = 0$). Then
\[
\E[f g] = \| f g \|_1 \geq \|f\|_p \|g\|_{p'},
\]
where equality holds if $g = f^{p/p'}$.
\end{lemma}

\begin{proof}[ of Corollary~\ref{cor:two-fcn-reverse-BB}]
By definition, the left-hand side of (\ref{eq:two_Borell}) is $\E[f T_\rho g]$.
We claim it suffices to prove (\ref{eq:two_Borell}) for
$\rho = (1-p)^{1/2}(1-q)^{1/2}$. Indeed, otherwise, let $r$ satisfy
$\rho = (1-p)^{1/2}(1-r)^{1/2}$ and note that $r \geq q$. Then,
assuming (\ref{eq:two_Borell}) holds for $p,r$ and $\rho$ we obtain:
\[
\E[f T_\rho g] \geq \|f\|_p \|g\|_r \geq \|f\|_p \|g\|_q,
\]
as needed.

We now assume $\rho = (1-p)^{1/2}(1-q)^{1/2}$.
Let $p'$ satisfy $1/p +1/p' = 1$.
Applying the reverse H\"{o}lder inequality we get that
$\E[f T_\rho g] \geq \|f\|_p \|T_\rho g\|_{p'}$.
Note that, since $1/(1-p') = 1-p$, the fact that $\rho =
(1-p)^{1/2}(1-q)^{1/2}$ implies $\rho = (1-q)^{1/2}(1-p')^{-1/2}$. Therefore, using
the reverse Bonami-Beckner inequality with $p' \leq q \leq 1$, we conclude that
\[
\E[f(x)g(y)] \geq \|f\|_p \|T_\rho g\|_{p'} \geq  \|f\|_p \|g\|_q.
\]
\end{proof}

\begin{lemma} \label{lem:duality}
It suffices to prove (\ref{eq:borell}) for $0 < q < p < 1$.
\end{lemma}

\begin{proof}
Note first that the case $p=1$ follows from the case $p < 1$ by continuity.
Recall that $1-p = \rho^2(1-q)$. Thus, $p > q$.
Suppose (\ref{eq:borell}) holds for $0 < q < p < 1$. Then by
continuity we obtain (\ref{eq:borell}) for $0 \leq q < p < 1$.
From $1-p = \rho^2(1-q)$, it follows that
$1-q' = 1/(1-q) = \rho^2/(1-p) = \rho^2(1-p')$. Therefore if $p \leq
0$, then $p'=1-1/(1-p) \geq 0$ and $q' = 1 - \rho^2/(1-p) > p'\geq
0$. We now conclude that if $f$ is non-negative, then
\begin{eqnarray*}
\|T_{\rho} f\|_q &=& \inf\{ \|g T_{\rho} f\|_1 : \|g\|_{q'} = 1, g \geq 0\}
\,\,\,\text{ (by reverse H\"{older}) } \\
&=& \inf\{ \|f T_{\rho} g\|_1 : \|g\|_{q'} = 1, g \geq 0\}
\,\,\,\text{ (by reversibility) } \\
&\geq& \inf\{ \|f\|_p \|T_{\rho} g\|_{p'} : \|g\|_{q'} = 1, g \geq 0\}
\,\,\,\mbox{ (by reverse H\"{older}) } \\
&\geq& \|f\|_p \inf\{ \|g\|_{q'} : \|g\|_{q'} = 1, g \geq 0\} = \|f\|_p
\,\,\,\mbox{ (by (\ref{eq:borell}) for } 0 \leq p' < q' < 1 \text{)}.
\end{eqnarray*}
We have thus obtained that (\ref{eq:borell}) holds for $p \leq 0$.
The remaining case is $p > 0 > q$.
Let $r=0$ and choose $\rho_1,\rho_2$ such that $(1-p) = \rho_2^2 (1-r)$
and $(1-r) = \rho_1^2 (1-q)$. Note that $0<\rho_1,\rho_2<1$ and
that $\rho = \rho_1 \rho_2$. The latter equality implies that
$T_{\rho} = T_{\rho_1} T_{\rho_2}$ (this is known as the ``semi-group
property''). Now
\[
\|T_{\rho} f\|_q = \|T_{\rho_1} T_{\rho_2} f\|_q  \geq \|T_{\rho_2} f\|_r
 \geq \|f\|_p,
\]
where the first inequality follows since $q < r \leq 0$ and the second
since $p > r \geq 0$.

We have thus completed the proof.
\end{proof}

\end{document}